\documentclass[12pt]{article}

\usepackage{amsfonts}
\usepackage{amsmath}
\usepackage{amsthm}
\usepackage{amssymb}
\usepackage{eucal}
\usepackage{mathtools}
\usepackage[colorlinks]{hyperref}
\usepackage{xcolor}
\numberwithin{equation}{section}

\renewcommand{\leq}{\leqslant}
\renewcommand{\geq}{\geqslant}
\newcommand{\Det}{\det}
\newcommand{\Tr}{\operatorname{Tr}}
\newcommand{\Area}{\operatorname{Area}}
\newcommand{\Length}{\operatorname{Length}}
\newtheorem{proposition}{Proposition}[section]
\newtheorem{lemma}[proposition]{Lemma}
\newtheorem{theorem}[proposition]{Theorem}

\newtheorem{remark}[proposition]{Remark}

\title{ Determinants of Dirichlet and Neumann Laplacians
on Constant Curvature Surfaces with Geodesic Polygonal Boundary}
\author{Victor Kalvin}

\date{}
\begin{document}
\maketitle

\begin{abstract}
We express the zeta regularized determinants of the Dirichlet and Neumann Laplacians on compact Riemann surfaces
 with constant curvature metrics and geodesic polygonal boundary in terms of the Laplacian determinant on their closed conical doubles.
  The formulae apply without restrictions on the genus or number of boundary components. They involve explicit corner angle factors and a quotient
of the Reidemeister and $L^2$ norms on the determinant line of
cohomology. The proof combines a new relative BFK type gluing formula with the comparison theorem for analytic and combinatorial torsion on manifolds with boundary, due to L\"uck and Vishik. As an application, we give a rigorous proof of the Aurell-Salomonson formula for Euclidean polygons. We also obtain explicit determinants for spherical and hyperbolic triangles and certain other geodesic polygons from known formulae for their doubles and discuss some extremal properties.
\end{abstract}

%\noindent\textbf{MSC:} 58J52.\par
%\noindent\textbf{Keywords:} zeta regularized determinants; relative determinants;
%analytic surgery; BFK gluing; conical singularities; analytic and Reidemeister torsion.\par

\section{Introduction and main result}
\label{Sec1}

%We consider compact Riemann surfaces with constant curvature metrics and geodesic polygonal boundary. We express the Dirichlet and Neumann determinants in terms of the  determinant on the closed conical double. Thus every
%determinant formula for the double
%yields the determinants of both the Dirichlet and Neumann Laplacians on
%the original surface.

%------------
We consider compact Riemann surfaces with constant curvature metrics and geodesic polygonal boundary. We express the  determinants of the Friedrichs Dirichlet and Neumann Laplacians in terms of the determinant of the Friedrichs Laplacian on the closed conical double.

The singular Polyakov and Polyakov--Alvarez anomaly formulae in~\cite{Ref11} led to closed explicit expressions for the determinants on constant curvature double triangles~\cite{Ref14} and triangulated conical spheres~\cite{Ref15}. For doubles of genus greater than one, there is a general determinant formula~\cite[Proposition 3.9]{Ref11}, with the previously undetermined universal constant evaluated in~\cite[Remark 7.2]{Ref17}. Thus determinant formulae without undetermined constants are available for all the doubles considered here, although not always in a closed explicit form. The obstacle to general closed explicit expressions is the longstanding problem of explicit uniformization,  which is resolved only in certain special cases, often possessing additional symmetries.

A BFK gluing formula~\cite{Ref5} would relate the Dirichlet determinants on the two halves to the determinant on the double. In the presence of boundary corners, however, the separating cut passes through conical points, and the usual BFK gluing formulae do not apply~\cite{Ref11,Ref15}.  Here we address this difficulty for reflection doubles, expressing the Dirichlet and Neumann determinants of each half in terms of the determinant of the double, with explicit corner angle and torsion factors. Thus every determinant formula for the double yields the determinants of both the Dirichlet and Neumann Laplacians on the original surface.

Aurell and Salomonson obtained explicit determinant formulae for polyhedral surfaces and Euclidean polygons using partially heuristic arguments~\cite{Ref3,Ref4}. Their formula for polyhedral surfaces of spherical topology was rigorously derived as a byproduct in~\cite[Section 3.2]{Ref11}; in this paper we deduce their polygonal formula as an application of our main theorem. 

The results for double triangles~\cite{Ref14} and regular dihedra, that is, doubles of regular polygons~\cite[Examples 5.2 and 5.7]{Ref15}, also lead to
%The results of~\cite{Ref14,Ref15} for the doubles also lead to 
closed explicit formulae for the determinants of Friedrichs Dirichlet and Neumann Laplacians on flat, spherical, and hyperbolic triangles and further geodesic polygons. For fixed curvature and area, the stationarity results for the corresponding doubles~\cite{Ref14,Ref15} imply that equilateral triangles and regular Euclidean and hyperbolic polygons are stationary points of both the Dirichlet and Neumann determinants. 
Moreover,  at fixed curvature and area, the equilateral triangle is a strict local minimum of both the Dirichlet and Neumann determinants if the area is sufficiently small. For larger areas, the situation can change. For example, among Euclidean triangles of unit area, the equilateral triangle is a strict local maximum of the Dirichlet determinant but a strict local minimum of the Neumann determinant. For Euclidean and hyperbolic triangles, both determinants are unbounded from above at fixed curvature and area, tending to infinity as an angle tends to zero. These extremal properties contrast sharply with the classical determinant maximization results for smooth metrics of Osgood, Phillips, and Sarnak~\cite{Ref31}. %Both attain absolute minima; in particular, among Euclidean triangles of unit area the Dirichlet determinant attaints its minimum on a non-equilateral triangle.

Some other determinant formulae for closed conical Riemann surfaces can be found in~\cite{Ref22,Ref18,Ref19,Ref10,Ref12,Ref17}, and for surfaces with boundary in~\cite{Ref2,Ref7,Ref11,Ref13}.

%---------------

%Determinant comparison formulae and explicit expressions are available for
%closed conical Riemann surfaces~\cite{Ref22,Ref18,Ref19,Ref10,Ref11,Ref14,Ref12,Ref15,Ref17} and for conical surfaces with smooth
%boundary~\cite{Ref11,Ref13}. %The constant in the principal-stratum formula~\cite{Ref22} was evaluated for all genera in~\cite[Remark~7.2]{Ref17}.
% For Polyakov--Alvarez comparison formulae on curvilinear polygons under conformal changes smooth up to the corners, see~\cite[Theorem~1.7]{Ref2} and~\cite[(8)]{Ref7}.  Aurell and Salomonson obtained explicit
%determinant formulae for polyhedral surfaces and Euclidean polygons using
%partially heuristic arguments~\cite{Ref3,Ref4}.
%Aurell-Salomonson formula for polyhedral surfaces of spherical topology was rigorously
%derived as a byproduct in~\cite[Section~3.2]{Ref11}; in this paper we deduce their polygonal
%formula as an application of our main theorem.
%
%

Let $M$ be a compact  Riemann surface of genus
$g$, with $\ell\geq1$ boundary components, and let $m$ be a
conformal metric on $M$. Assume that the metric $m$ is
smooth and has constant curvature in the interior, and that the
boundary is geodesic except at finitely many vertices $p_j$ with
interior angles $\pi\alpha_j>0$.  We call $(M,m)$ a constant curvature surface with geodesic polygonal boundary.
The assumption that $m$ is smooth in the interior of $M$ is made only
for simplicity; finitely many interior conical singularities can be
allowed without additional difficulty.

The double
$$
 X=M\cup_{\partial M}M
$$
is obtained by identifying corresponding boundary points of
$M$ and a copy endowed with the conjugate complex structure.
The closed Riemann surface $X$ has genus $2g+\ell-1$.
The metric $m$ extends by reflection to a constant curvature
conformal metric on $X$, with conical angles $2\pi\alpha_j$
at $p_j\in\partial M$. The isometric antiholomorphic
reflection $\iota:(X,m)\to(X,m)$ exchanges the two copies of  $(M,m)$ and
fixes their common boundary.
The corners of $\partial M$ and the conical singularities of $X$ are encoded entirely in the metric $m$.

Let $\Delta_m$ be the  Friedrichs selfadjoint
extension in $L^2(X,m)$ of the  Laplace--Beltrami operator
on the functions in  $C_c^\infty(X\setminus\bigcup_j\{p_j\})$. Similarly, let $\Delta_m^D$ and $\Delta_m^N$ be the Friedrichs selfadjoint
extensions of the Laplace--Beltrami operator on $M$, with Dirichlet
and Neumann boundary conditions, respectively.

All three Laplacians have discrete spectrum. We define the spectral zeta regularized determinant of $\Delta_m$ in the standard way:
$$
 \zeta_m(s):=\sum_{\substack{\lambda_k\in\operatorname{spec}(\Delta_m), \lambda_k\neq 0}}
 \lambda_k^{-s},\quad \Re s>1,
 \qquad
 \Det\Delta_m:=\exp\left(
 -\left.\partial_s\zeta_m(s)\right|_{s=0}\right),
$$
where the eigenvalues are counted with multiplicities and the
derivative is taken after analytic continuation. The  determinants of
$\Delta_m^D$ and $\Delta_m^N$ are defined in the same way. The Laplacian  $\Delta_m^D$ is strictly
positive, whereas the kernels of $\Delta_m$ and $\Delta_m^N$
consist of constant functions.

Let $H^q(M)$ denote the ordinary de Rham cohomology
with complex coefficients: closed smooth $q$-forms modulo exact
forms, without boundary conditions.  Since $M$ is connected and has nonempty
boundary,
$$
 H^0(M)=\mathbb C,\qquad
 \dim H^1(M)=2g+\ell-1,\qquad
 H^2(M)=0.
$$
In the convention of \cite[Sec.~1.1]{Ref30}, the determinant line
of this cohomology is
$$
 \det H^\bullet(M)
 =\bigl(\det H^0(M)\bigr)^{-1}
  \otimes\det H^1(M).
$$
Here $ \det H^q=\bigwedge^{\dim H^q}H^q$ stands for the one-dimensional vector space, 
and the exponent $-1$ denotes its
dual.  Equip $H^0(M)$ with the $L^2(M,m)$-norm of constant functions,
and $H^1(M)$ with the $L^2(M,m)$-norm of its absolute
harmonic representatives in the conformal structure of
$M$. These are the one-forms $\omega$ satisfying
$d\omega=d(*\omega)=0$ whose $*\omega$ has zero pullback to
$\partial M$; the Hodge star on one-forms is determined by the
conformal structure. The resulting norm on the determinant line $\det H^\bullet(M)$
is denoted by $\|\cdot\|_{\det H^\bullet(M)}$.

Let $\tau_0(M)$ be the
squared Reidemeister (combinatorial) torsion norm on $\det H^\bullet(M)$,
with complex coefficients; see~\cite[(1.6)]{Ref30}. It is obtained from
the cellular cochain complex of a triangulation, with the simplex
cochains forming an orthonormal basis, and is independent of the triangulation.

 For any nonzero element $v$ of
the line   $\det H^\bullet(M)$, the quotient $ {\|v\|_{\tau_0(M)}^2}/
{ \|v\|_{\det H^\bullet(M)}^2}$
is $v$-independent.  If, for example,  $M$ is a topological disk, then the quotient reduces to the area of the geodesic polygon $(M,m)$. 

\begin{theorem}\label{Thm1.1}
Under the assumptions above, the spectral zeta regularized determinants of
the Friedrichs  Dirichlet and Neumann Laplacians on the constant curvature surface  $(M,m)$ with geodesic polygonal boundary $\partial M$ are given by
\begin{equation}\label{eq1.1}
 \begin{aligned}
 \Det\Delta_m^D
 &=\left(
 \frac{\|v\|_{\det H^\bullet(M)}^2}
 {\|v\|_{\tau_0(M)}^2}\,\Det\Delta_m
 \right)^{1/2}\prod_j\alpha_j^{1/4},\\[4pt]
 \Det\Delta_m^N
 &=\left(
 \frac{\|v\|_{\tau_0(M)}^2}
 {\|v\|_{\det H^\bullet(M)}^2}\,\Det\Delta_m
 \right)^{1/2}\prod_j\alpha_j^{-1/4}.
 \end{aligned}
\end{equation} 
Here $\Delta_m$ is the Friedrichs Laplacian on the closed Riemann surface $X=M\cup_{\partial M}M$, equipped with the doubled conical metric $m$. The interior angles at the corners of $\partial M$ are $\pi\alpha_j$.
 The element
$v$ of the determinant line $\det H^\bullet(M)$ is arbitrary and nonzero; the norms
$\|\cdot\|_{\tau_0(M)}$ and
$\|\cdot\|_{\det H^\bullet(M)}$ are the Reidemeister torsion
norm and the norm induced by $L^2(M,m)$, respectively.
\end{theorem}

We now outline the main ideas of the proof.
The decomposition into odd and even functions with respect to the reflection $\iota$ shows that the determinant on the double $\det\Delta_m$ is equal to the product of the Dirichlet $\det\Delta_m^D$ and Neumann $\det\Delta_m^N$ determinants. To recover them separately, it therefore suffices
to determine the ratio ${\Det\Delta_m^N}/{\Det\Delta_m^D}
$. 

For smooth metrics that are products near the boundary, the ratio follows from the extension of the Cheeger--M\"uller theorem to manifolds with boundary, due to L\"uck and Vishik~\cite{Ref26,Ref30}. The Polyakov and Polyakov--Alvarez formulae~\cite{Ref11} then extend the result  to reflection invariant metrics  that are smooth near the boundary. 

Corners at the boundary $\partial M$ require a different argument. We cut out small geodesic disks centred at the conical points of the double $X$ and smooth the metric inside them, preserving rotational symmetry. The BFK gluing formula~\cite{Ref5} could then compare the original metric with its smoothing~\cite{Ref11}, but the comparison is needed separately for odd and even functions. The usual approach to the proof of the BFK formula fixes the gluing constants from the large spectral parameter asymptotics of the individual Laplacian and Neumann jump determinants. Obtaining the required asymptotics for the odd and even restrictions seems to be problematic.

We instead develop a relative gluing approach and establish a new BFK type formula comparing the singular metric directly with its smoothing, separately for odd and even functions. This formula involves a relative Neumann jump operator that differs from the identity by a trace class operator. It therefore has a Fredholm determinant. Trace norm estimates, which pass directly to the odd and even restrictions, show that the corresponding Fredholm determinants tend to one as the spectral parameter goes to infinity. This fixes the relative gluing constants without the individual determinant asymptotics.

On the double with the disks removed, the metric is unchanged.  As a result, the corresponding contributions in the relative BFK formula cancel. Only the Dirichlet or Neumann determinants for the singular and smoothed metrics, their disk counterparts, and the Fredholm determinant of the relative Neumann jump remain.

We then pass to zero spectral parameter, treating the odd and even restrictions separately. The Fredholm determinant tends to one in the odd case, whereas in the even case the constant zero mode produces a ratio of areas.

To compare the disk determinants, we separate variables in the (rotationally symmetric) metric disks. For each metric, the contributions of nonzero angular frequencies cancel in the ratio of the even and odd disk determinants, since the corresponding radial spectra coincide. Comparing the remaining radial determinants for the conical and smoothed metrics  by using Lesch's formula~\cite{Ref25}, we obtain the explicit corner angle factors. 
Together with the smooth torsion formula, these comparisons determine the required ratio  ${\Det\Delta_m^N}/{\Det\Delta_m^D}
$, and we recover both determinants from their product.

The paper is organized as follows. In Section~\ref{Sec2} we construct a local smoothing of the conical metric and compare the corresponding disk determinants. Section~\ref{Sec3}
studies the relative Neumann jump operator and establishes the estimates
used to prove the relative BFK formula in Section~\ref{Sec4}. In Section~\ref{Sec5}  we pass to zero spectral parameter and complete the proof of Theorem~\ref{Thm1.1} using the smooth torsion calibration.
Section~\ref{Sec6} discusses BFK gluing and the torsion factor.
Section~\ref{Sec7} contains examples and applications.

%-----------------------------------

\section{Local determinant ratio under conical smoothing}\label{Sec2}

Let $p_j\in X$ be one of the conical points of the double, of full
angle $2\pi\alpha_j$.  In this section we fix $j$ and suppress the index in notation
internal to this section.  

By~\cite[Theorem~4.1]{Ref29}, in a small neighbourhood 
there exist geodesic polar coordinates
$$
 0<r<\epsilon,
 \qquad
 \theta\in\mathbb R/(2\pi\mathbb Z),
$$
centred at $p_j$, such that the metric $m$ takes the form
\begin{equation*}
 m=dr^2+\alpha_j^2\mathsf s_K^2(r,\theta)\,d\theta^2.
\end{equation*}
Here $\mathsf s_K\in C^\infty((0,\epsilon)\times
\mathbb R/(2\pi\mathbb Z))$ is the positive normalized angular coefficient in
Troyanov's polar parametrization.  Moreover, $\mathsf s_K$ is continuous at the  vertex and satisfies there the one-sided initial data
\begin{equation}\label{eq2.1}
 \mathsf s_K(0,\theta)=0,
 \qquad
 \partial_r^+\mathsf s_K(0,\theta)=1.
\end{equation}
 Since
$\mathsf s_K>0$ for $r>0$, the Gaussian curvature equation $K=-\partial_r^2\mathsf s_K/\mathsf s_K$ is equivalent to the scalar
Jacobi equation
$$
 \partial_r^2\mathsf s_K(r,\theta)
 +K\mathsf s_K(r,\theta)=0,
 \qquad 0<r<\epsilon,\qquad \theta\in\mathbb R/(2\pi\mathbb Z).
$$
The initial data \eqref{eq2.1} together with  unique solvability of the ordinary
differential equation imply that  $\mathsf s_K$ is
 $\theta$-independent and  given by$$
 \mathsf s_K(r)=
 \begin{cases}
  K^{-1/2}\sin(\sqrt K\,r),&K>0,\\[2pt]
  r,&K=0,\\[2pt]
  (-K)^{-1/2}\sinh(\sqrt{-K}\,r),&K<0.
 \end{cases}
$$
 
 Recall that on the double, the reflection $\iota$ exchanges the two copies of $M$ across their common boundary. Let the polar coordinates be chosen so that locally
$\iota(r,\theta)=(r,-\theta)$.

Choose and fix $\epsilon>0$, the same for all conical points, so
that the corresponding closed polar disks $\mathcal U^j$ are pairwise
disjoint and contain no other conical points.  For the fixed point
write $\mathcal U^j=\{r\leq\epsilon\}$.  The quotient of the conical metric disk $(\mathcal U^j,m)$ by
the reflection $\iota$ is the geodesic sector
$$
 \mathcal Q^j=\mathcal U^j/\iota=\{0\leq r\leq \epsilon,\ 0\leq\theta\leq\pi\}.
$$
The boundary of $\mathcal Q^j$ consists of  the arc $\gamma=\{r=\epsilon, \theta\in[0,\pi] \}$ of the geodesic circle  centered at $p_j$, and 
 the two radial geodesic parts $\sigma$. At  the vertex $p_j$, corresponding to $r=0$,  the metric sector $(\mathcal Q^j, m)$ has a corner of angle $\pi\alpha_j$.

Now we are about to smooth $m$ near $p_j$, preserving the rotational and
reflection symmetries and leaving the metric unchanged near the
outer arc $\gamma$.  The resulting smooth  conformal metric will be denoted by $\widetilde m$.

We  first define the conformal radius 
$$
 \rho(r)
 =
 \exp\left(
 -\int_r^\epsilon
 \frac{dt}{\alpha_j\mathsf s_K(t)}
 \right),
$$
so that $\rho(\epsilon)=1$. 
 Then $z=\rho e^{i\theta}$ is a local holomorphic parameter centered at $p_j$,
$\iota(z)=\bar z$, and
$$
 m
 =
 \rho^{2\alpha_j-2}e^{2\phi(\rho)}|dz|^2,
 \qquad
 \phi(\rho)
 =
 \log\frac{\alpha_j\mathsf s_K(r(\rho))}{\rho^{\alpha_j}}.
$$
Choose a cutoff function $\eta\in C^\infty([0,1])$ such that
$\eta=0$ near $0$ and $\eta=1$ near $1$. Introduce the smoothed metric 
$$
 \widetilde m
 =
 \exp\left(
 2\eta(\rho)
 \bigl((\alpha_j-1)\log\rho+\phi(\rho)\bigr)
 \right)|dz|^2.
$$
Thus $\widetilde m=|dz|^2$ near $p_j$ and
$\widetilde m=m$ near $\gamma$.  Moreover, on
$\mathcal U^j\setminus\{p_j\}$, we have
$$
 \widetilde m=e^{2\omega}m,
 \qquad
 \omega
 =
 (\eta-1)
 \bigl((\alpha_j-1)\log\rho+\phi(\rho)\bigr).
$$

For $h\in\{m,\widetilde m\}$, let
$\Delta_{\mathcal U^j,h}$ be the Friedrichs Dirichlet Laplacian on the metric
disk $(\mathcal U^j,h)$.  The reflection $\iota$ induces the unitary involution
$$
 \iota^*:L^2(\mathcal U^j,h)\longrightarrow L^2(\mathcal U^j,h),
 \qquad
 (\iota^*u)(z)=u(\bar z).
$$
The operator
$\Delta_{\mathcal U^j,h}$ commutes with $\iota^*$ as the quadratic form and its domain
are invariant under $\iota^*$.  Denote 
$$
 \Pi^D=\frac12(I-\iota^*),
 \qquad
 \Pi^N=\frac12(I+\iota^*),
$$
and define
$$
 \Delta_{\mathcal U^j,h}^D
 :=
 \left.\Delta_{\mathcal U^j,h}\right|_{\operatorname{Ran}\Pi^D},
 \qquad
 \Delta_{\mathcal U^j,h}^N
 :=
 \left.\Delta_{\mathcal U^j,h}\right|_{\operatorname{Ran}\Pi^N},
$$
where $\operatorname{Ran}\Pi^D$ and
$\operatorname{Ran}\Pi^N$ stand for  the ranges of the two orthogonal projections $\Pi^D$ and $\Pi^N$ in $L^2(\mathcal U^j,h)$.
These restrictions are unitarily equivalent to the Laplacians on
$(\mathcal Q^j,h)$ with, respectively, Dirichlet and Neumann
conditions on $\sigma$, and with the Dirichlet condition on
$\gamma$ in both cases.

\begin{proposition}\label{Prop2.1}  
For $h\in\{m,\widetilde m\}$ and
$\sharp\in\{D,N\}$, the spectral zeta function
$$
 \zeta_{\mathcal U^j,h}^\sharp(s)
 :=
 \Tr\bigl(\Delta_{\mathcal U^j,h}^\sharp\bigr)^{-s},
 \qquad \Re s>1,
$$
admits an analytic continuation to a neighborhood of $s=0$. Let 
$$
 \Det\Delta_{\mathcal U^j,h}^\sharp
 :=
 \exp\left(
 -\left.\partial_s
 \zeta_{\mathcal U^j,h}^\sharp(s)\right|_{s=0}
 \right)
$$
stand for the zeta regularized spectral  determinant. 
Denote 
 \begin{equation*}
 F_j(h):=
 \frac{\Det\Delta_{\mathcal U^j,h}^N}
      {\Det\Delta_{\mathcal U^j,h}^D},\quad h\in\{m,\widetilde m\}. 
\end{equation*}
Then the conical metric $m$ and its rotational smoothing $\widetilde m$ satisfy
\begin{equation}\label{eq2.2}
 \frac{F_j(m)}{F_j(\widetilde m)}
 =\alpha_j^{-1/2}.
\end{equation}
\end{proposition}

\begin{proof} 
Since both metrics are rotationally invariant,   the  boundary $\partial\mathcal U^j=\{\rho=1\}$ is a centered
geodesic circle for either metric.  Fix
$h\in\{m,\widetilde m\}$, let $r$ denote the geodesic distance
from $p_j$ in the metric $h$, and denote the  value of $r$  on
$\partial\mathcal U^j$ by $\epsilon$.  For $h=m$, these
are the coordinate $r$ and the number $\epsilon$ fixed above;
for $h=\widetilde m$, the same symbols denote the corresponding
geodesic quantities for the smoothed metric.  

The rotationally symmetric metric $h$ takes the form
\begin{equation*}
 h=dr^2+f(r)^2\,d\theta^2,
 \qquad
 0<r<\epsilon,
 \qquad
 \theta\in\mathbb R/2\pi\mathbb Z,
\end{equation*}
At the centre, the function $f$ has the asymptotics
\begin{equation*}
 f(r)=\beta r+O(r^3),
 \qquad r\to0,
 \qquad
 \beta=
 \begin{cases}
  \alpha_j,&h=m,\\
  1,&h=\widetilde m.
 \end{cases}
\end{equation*}
More precisely,  we have 
\begin{equation}\label{eq2.3}
 f(r)=\beta r\bigl(1+r^2a(r)\bigr),
 \qquad
 a\in C^\infty[0,\epsilon],
 \qquad
 f(r)>0\quad(0<r\leq\epsilon),
\end{equation}
\begin{equation}\label{eq2.4}
 2\pi f(\epsilon)=\Length(\partial\mathcal U^j,h).
\end{equation}
Since $m=\widetilde m$ near $\partial\mathcal U^j$, one also has
$$
 \Length(\partial\mathcal U^j,m)
 =
 \Length(\partial\mathcal U^j,\widetilde m).
$$

We first show that the spectral zeta function  $\zeta_{\mathcal U^j,h}^\sharp$   admits an analytic continuation to $s=0$. Indeed, the
spectrum of  the Dirichlet Laplacian $\Delta_{\mathcal U^j,h}$ is discrete and its spectral zeta function $\zeta_{\mathcal U^j,h}$
is regular at zero, e.g.~\cite[Sec.~1.1]{Ref11}.  In the space 
$L^2((0,\epsilon),f(r)\,dr)$ consider the selfadjoint Friedrichs extension $L_k$ of the operator
$$
 -\frac1{f(r)}\frac d{dr}f(r)\frac d{dr}
       +\frac{k^2}{f(r)^2}, \quad k=0,1,2,\dots, 
$$
with the homogeneous Dirichlet condition
at $r=\epsilon$.  Under the reflection
$\iota(r,\theta)=(r,-\theta)$, the circular Fourier modes split into
the odd modes $\sin(k\theta)$, $k\geq1$, and the even modes
$1,\cos(k\theta)$.  Their restrictions to the sector $\mathcal Q^j$ satisfy,
respectively, the Dirichlet and Neumann conditions on $\sigma$. The radial operator in the $k$-th mode is $L_k$.
With respect to the corresponding unitary Fourier decomposition of
$L^2(\mathcal U^j,h)$, the Laplacians are unitarily
equivalent to the following orthogonal sums of selfadjoint operators:
$$
 \Delta_{\mathcal U^j,h}
 \simeq L_0\oplus\bigoplus_{k\geq1}(L_k\oplus L_k),
$$
whereas its restrictions are
$$
 \Delta_{\mathcal U^j,h}^D\simeq\bigoplus_{k\geq1}L_k,
 \qquad
 \Delta_{\mathcal U^j,h}^N\simeq L_0\oplus\bigoplus_{k\geq1}L_k.
$$
The spectral series for $\zeta_{\mathcal U^j,h}$ converges absolutely for
$\Re s>1$; hence so do its odd and even subseries.  The orthogonal
decompositions above therefore imply
\begin{equation}\label{eq2.5}
 \zeta_{\mathcal U^j,h}^D(s)
 =\frac12\bigl(\zeta_{\mathcal U^j,h}(s)-\zeta_{0,h}(s)\bigr),
 \qquad
 \zeta_{\mathcal U^j,h}^N(s)
 =\frac12\bigl(\zeta_{\mathcal U^j,h}(s)+\zeta_{0,h}(s)\bigr),
\end{equation}
where $\zeta_{0,h}$ is the spectral zeta function of $L_0$.
To pass to $s=0$, use the unitary map
$$
 u\mapsto y=\sqrt f\,u:
 L^2((0,\epsilon),f(r)\,dr)\longrightarrow L^2(0,\epsilon)
$$
which transforms $L_0$ into
$
 H_f=-\partial_r^2+{(\sqrt f)''}/{\sqrt f}$. By \eqref{eq2.3} we have 
\begin{equation}\label{eq2.6}
 \frac{(\sqrt f)''}{\sqrt f}
 =-\frac1{4r^2}+V_f(r),
 \qquad
 V_f\in C^\infty[0,\epsilon].
\end{equation}
Thus $H_f$ is regular singular  in the sense of Lesch~\cite{Ref25}. By~\cite[Theorem~1.1]{Ref25} the zeta function $\zeta_{0,h}$ admits an analytic
continuation to a neighborhood of $s=0$. 

 Since
\eqref{eq2.5} holds for $\Re s>1$ and  $\zeta_{\mathcal U^j,h}$ also
admits an analytic continuation to a neighborhood of $s=0$, the 
right hand sides of~\eqref{eq2.5} give the analytic continuations of
$\zeta_{\mathcal U^j,h}^D$ and $\zeta_{\mathcal U^j,h}^N$ to the same neighborhood.  Hence
$\Det\Delta_{\mathcal U^j,h}^D$ and $\Det\Delta_{\mathcal U^j,h}^N$ are well-defined.
Differentiating \eqref{eq2.5} at $s=0$, we obtain
\begin{equation*}
 F_j(h)=\Det L_0.
\end{equation*}

To evaluate $\Det L_0$, we use Lesch's determinant formula
\cite[Theorem~1.2]{Ref25}.  In the present case, the solution entering
this formula is the unique function $\phi_f$ satisfying
$$
 H_f\phi_f=0,
 \qquad
 \phi_f(r)=\sqrt r+o(\sqrt r)
 \quad\text{as }r\to0.
$$
Since $H_f\sqrt f=0$ and
$f(r)=\beta r(1+O(r^2))$, this solution is
\begin{equation}\label{eq2.7}
 \phi_f(r)=\sqrt{{f(r)}/{\beta}}.
\end{equation}

In the notation of Lesch, the coefficient $-1/4$ in~\eqref{eq2.6}
 corresponds to $\nu_0=0$, whereas the Dirichlet
condition at the regular endpoint $r=\epsilon$ corresponds to
$\nu_1=1/2$; see 
\cite[(1.26) and (1.35)]{Ref25}.  Let $\psi_f$ be the solution of
$H_f\psi_f=0$ normalized at $r=\epsilon$ as in
\cite[(1.38)]{Ref25}.  Since $\nu_1=1/2$, this normalization reads
$$
 \psi_f(\epsilon)=0,
 \qquad
 \psi_f'(\epsilon)=-1.
$$
The Wronskian in Lesch's formula can be evaluated at
$r=\epsilon$. We have
$$
 W(\psi_f,\phi_f)
 =\psi_f(\epsilon)\phi_f'(\epsilon)
  -\psi_f'(\epsilon)\phi_f(\epsilon)
 =\phi_f(\epsilon).
$$
As a result, Theorem~1.2 of \cite{Ref25} implies
\begin{equation}\label{eq2.8}
 \Det H_f
 =
 \frac{\pi W(\psi_f,\phi_f)}
 {2^{\nu_0+\nu_1}
  \Gamma(\nu_0+1)\Gamma(\nu_1+1)}
 =
 \frac{\pi\phi_f(\epsilon)}
 {2^{1/2}\Gamma(1)\Gamma(3/2)}
 =
 \sqrt{2\pi}\,\phi_f(\epsilon).
\end{equation}
Although Theorem~1.2 of \cite{Ref25} is stated on $[0,1]$, its
result is unchanged under the linear rescaling $r=\epsilon x$.
Indeed, for the rescaled operator,
$\zeta(0)=-(\nu_0+\nu_1)/2=-1/4$ by
\cite[Theorem~1.1]{Ref25}, so the determinant acquires the factor
$\epsilon^{1/2}$, which cancels the factor $\epsilon^{-1/2}$ in
the normalized solution.

 Substitution of
\eqref{eq2.7} into~\eqref{eq2.8},
together with~\eqref{eq2.4}, now yields
\begin{equation*}
 F_j(h)=\sqrt{\frac{\Length(\partial\mathcal U^j,h)}{\beta}}.
\end{equation*}
For $h=m$, one has $\beta=\alpha_j$, whereas for
$h=\widetilde m$, one has $\beta=1$. Since the two boundary lengths are equal, this proves~\eqref{eq2.2}.
\end{proof}

\begin{remark}\label{Rem2.2}
For a flat conical disk of radius $\epsilon$, the identities
\eqref{eq2.5} and
\eqref{eq2.8} imply
$$
 \Det\Delta_{\mathcal U^j,m}
 =\sqrt{2\pi\epsilon}\,
   \bigl(\Det\Delta_{\mathcal U^j,m}^D\bigr)^2.
$$
This recovers the relation between the Dirichlet determinants on
a finite Euclidean sector and its conical double derived 
in~\cite[Sec.~3]{Ref2}; see also~\cite{Ref7}.
Our calculation also applies to the conical metrics admissible in the sense
of~\cite[Definition~2.1]{Ref11}, provided that they are rotationally
symmetric and their radial coefficient $f$
satisfies~\eqref{eq2.3}; constant curvature is not required.
\end{remark}

%++++++++++++++++++++++++++++++++++++++++++++++++
\section{Relative Neumann jump operator}\label{Sec3}

Let $\mathcal U^j$ be the closed disk of radius $\epsilon$ centered at
the conical point $p_j$ of the double. We assume that the radius $\epsilon$ was
chosen so that the disks $\mathcal U^j$ are pairwise disjoint.  In this
and the subsequent sections, we use the notation
$$\mathcal U:=\bigcup_j\mathcal U^j,\quad 
 X_\epsilon
: =X\setminus\operatorname{Int}\mathcal U.
$$
Let $\partial X_\epsilon=\partial \mathcal U$ stand for the boundary of
$X_\epsilon$.  On each $\mathcal U^j$, replace the metric $m$
by its smoothing $\widetilde m$ as in the preceding section, and then set
$\widetilde m=m$ on $X_\epsilon$.  Since the smoothed
metric $\widetilde m$  coincides with $m$ near every  cut $\partial\mathcal U^j$, this defines a 
smooth conformal reflection invariant metric $\widetilde m$ on $X$.  

As noted in the introduction, we assume only for simplicity that there are no conical points in  $\operatorname{Int} M$; to include them, one would  smooth the  boundary singularities of $m$, leaving the interior singularities unchanged. 

 Let
$u_h(\lambda)$, $\lambda>0$, be the unique weak solution of the Dirichlet boundary value problem 
\begin{equation}\label{eq3.1}
 (\Delta_h+\lambda)u_h(\lambda)=0
 \quad\text{on}\quad
 X\setminus
 \partial X_\epsilon,
 \qquad
 u_h(\lambda)=\varphi
 \quad\text{on}\quad
 \partial X_\epsilon,
\end{equation}
where $\varphi\in C^\infty(\partial X_\epsilon)$ and  $h$ stands for either of the two metrics  $m$ or $\widetilde m$.

Let $\nu$ be the outward unit normal to
$(\mathcal U,m)$, and let $ds$ denote the arclength measure on
$\partial X_\epsilon$.
Since the metrics coincide near $\partial X_\epsilon$, both $\nu$ and $ds$
are common to the metrics $m$ and $\widetilde m$.  All function spaces on
$\partial X_\epsilon$ are defined with respect to the common induced metric.
Define the Neumann jump operator
$R_h(\lambda)$ by the equality 
$$
 R_h(\lambda)\varphi
 =
 \partial_\nu
 \bigl(u_h(\lambda)|_{\mathcal U}\bigr)
 -
 \partial_\nu
 \bigl(u_h(\lambda)|_{X_\epsilon}\bigr).
$$
For $\lambda>0$, the operator
$$
 R_h(\lambda):
 C^\infty(\partial X_\epsilon)
 \rightarrow
 C^\infty(\partial X_\epsilon)
$$
is an invertible first-order elliptic classical pseudodifferential
operator.  We view both operators  $R_m(\lambda)$ and $R_{\widetilde m}(\lambda)$ as unbounded operators in
$L^2(\partial X_\epsilon)$, with domain
$H^1(\partial X_\epsilon)$.

The reflection $\iota$ acts unitarily on the functions in $L^2(\partial X_\epsilon)$ by
$(\iota^*\varphi)(x)=\varphi(\iota(x))$.  It is an isometry for
either metric and preserves each side of the cut.  Thus
$\iota^*u_h(\lambda)$ is the unique  solution to  the boundary value problem~\eqref{eq3.1} 
with $\varphi$ replaced by $\iota^*\varphi$.  Since $\iota$ preserves the
normal, every function $v$ on either side of the cut satisfies
$$
 \partial_\nu(\iota^*v)(x)
 = (\iota^*\partial_\nu v)(x).
$$
Consequently, for the Neumann jump operator  we have
$$
 R_h(\lambda)(\iota^*\varphi)
 =
\iota^* (R_h(\lambda)\varphi).
$$
 Thus the operator $R_h(\lambda)$ commutes with $\iota^*$ and therefore also with the orthogonal projections
 $$
 \Pi^D=\frac12(I-\iota^*),
 \qquad
 \Pi^N=\frac12(I+\iota^*)
$$
 in  $L^2(\partial X_\epsilon)$. For $\sharp\in\{D,N\}$, denote  by  $R_h^\sharp(\lambda)$ the restriction
of $R_h(\lambda)$ to
$\Pi^\sharp L^2(\partial X_\epsilon)$, with domain
$\Pi^\sharp H^1(\partial X_\epsilon)$.

Define the relative Neumann jump operator in  $L^2(\partial X_\epsilon)$ via
\begin{equation}\label{eq3.2}
 T(\lambda):=
 R_{\widetilde m}(\lambda)R_m(\lambda)^{-1}
 =
 I+\bigl(R_{\widetilde m}(\lambda)-R_m(\lambda)\bigr)
 R_m(\lambda)^{-1}.
\end{equation}
Since $R_m(\lambda)$ and $R_{\widetilde m}(\lambda)$ commute with orthogonal projections
$\Pi^D$ and $\Pi^N$, we have
\begin{equation}\label{eq3.3}
 T(\lambda)=T^D(\lambda)\oplus T^N(\lambda);\qquad
 T^\sharp(\lambda)=
 R_{\widetilde m}^\sharp(\lambda)
 \bigl(R_m^\sharp(\lambda)\bigr)^{-1},\quad \sharp\in\{D,N\}.
\end{equation}
%------

\begin{lemma}\label{Lem3.1} We write  $\mathcal B(E,F)$ for the space of bounded linear
operators from $E$ to $F$, equipped with the operator norm,
and $\mathfrak S_1$ for the trace class, equipped with
the trace norm.
Let $ S_\delta$ with some $\delta>0$ stand for the sector
$$
 S_\delta=\bigl\{\lambda\in\mathbb C:
 \Re\lambda>0,\quad
 |\Im\lambda|<\delta\Re\lambda\bigr\}.
$$
Then the following assertions hold true.
\begin{enumerate}
\item[1.] For $h\in\{m,\widetilde m\}$, the following
operator families are analytic:
\begin{equation}\label{eq3.4}
 S_\delta\ni\lambda\longmapsto R_h(\lambda)
 \in\mathcal B\bigl(
 H^1(\partial X_\epsilon),
 L^2(\partial X_\epsilon)
 \bigr),
\end{equation}
\begin{equation}\label{eq3.5}
 S_\delta\ni\lambda\longmapsto
 \bigl(R_{\widetilde m}(\lambda)-R_m(\lambda)\bigr)
 \in\mathcal B\bigl(
 H^1(\partial X_\epsilon),
 H^2(\partial X_\epsilon)
 \bigr).
\end{equation}
Moreover, both families also extend analytically to a neighborhood of $\lambda=0$.

\item[2.] The operator family
\begin{equation}\label{eq3.6}
S_\delta\ni \lambda\longmapsto T(\lambda)-I
 \in\mathfrak S_1\bigl(
 L^2(\partial X_\epsilon)\bigr)
\end{equation}
is analytic. Moreover, there exists $\lambda_0\geq1$ such that, for every
$n\in\mathbb N_0$, there is $C_n>0$ satisfying
\begin{equation}\label{eq3.7}
 \|T(\lambda)-I\|_{\mathfrak S_1}
 \leq C_n|\lambda|^{-n},
 \qquad
 \lambda\in S_\delta,\quad |\lambda|\geq\lambda_0.
\end{equation}

\item[3.] The previous  assertion also holds for $T(\lambda)$
replaced with its restrictions $T^\sharp(\lambda)$, 
 $\sharp\in\{D,N\}$. Moreover, the Fredholm determinant
$\det\nolimits_F T^\sharp(\lambda)$ does not vanish
on $S_\delta$, and
\begin{equation}\label{eq3.8}
 \det\nolimits_F T^\sharp(\lambda)\to 1
 \quad\text{as }|\lambda|\to\infty,\quad
 \lambda\in S_\delta.
\end{equation}
For the  branch of its logarithm tending to zero
as $|\lambda|\to \infty$, we have
\begin{equation}\label{eq3.9}
 \bigl|\log\det\nolimits_F T^\sharp(\lambda)\bigr|
 \leq C_n|\lambda|^{-n},
 \qquad
 \lambda\in S_\delta,\quad |\lambda|\geq\lambda_0,
 \quad n\in\mathbb N_0.
\end{equation}
\end{enumerate}
\end{lemma}
%------

\begin{proof}
The  metrics $m$  and $\widetilde m$ coincide in a fixed collar of
$\partial X_\epsilon$.  On each disk $\mathcal U^j$, let $z$
be the local holomorphic parameter introduced in the preceding section.   We shall use two smooth cutoff functions $\chi$ and $\eta$, such that
$\chi=\chi(|z|)$ and $\eta=\eta(|z|)$ on each disk $\mathcal U^j$. Recall that
$|z|=1$ on $\partial\mathcal U^j$.  Choose  $\chi\in C^\infty([0,1])$ such that $\chi=0$ whenever $m\neq\widetilde m$ and $\chi=1$ near $|z|=1$. The 
commutator $[\Delta_m,\chi]=[\Delta_{\widetilde m},\chi]$  is a first order differential operator with smooth coefficients supported
away from the boundary $\partial X_\epsilon$ and the
conical points.
 Choose $\eta$ satisfying $0\leq \eta\leq 1$,  $\eta=1$
near the support of $[\Delta_m,\chi]$, and $\eta=0$ near $z=0$
and $|z|=1$.  We extend $\chi$ by one and $\eta$ by zero to $X_\epsilon$. Clearly, $\eta[\Delta_m,\chi] \eta=[\Delta_m,\chi]$.   Below all function spaces and norms on $\mathcal U$ and its subdomains are defined with respect to the smooth metric
$\widetilde m$.

We first prove the analyticity of $R_h(\lambda)$ in the sector $S_\delta$ by establishing the analyticity of the difference $R_h(\lambda)-R_h(0)$.

 Let
$\Delta_{\partial X_\epsilon,h}$ with $h\in\{m,\widetilde m\}$ be the Friedrichs
Dirichlet Laplacian on
$X\setminus\partial X_\epsilon$. For $\varphi\in H^1(\partial X_\epsilon)$, the solutions  $u_h(\lambda)$ and $u_h(0)$ of the Dirichlet problem~\eqref{eq3.1} satisfy
\begin{equation}\label{eq3.10}
u_h(\lambda)-u_h(0)
=-\lambda
(\Delta_{\partial X_\epsilon,h}+\lambda)^{-1}
u_h(0).
\end{equation}
Since the harmonic extension $\varphi\mapsto u_h(0)$ is bounded into $L^2(X,h)$, the right hand side of~\eqref{eq3.10} is analytic on $S_\delta$ as a family of bounded operators from $H^1(\partial X_\epsilon)$ to $L^2(X,h)$.
 The standard local elliptic a priori estimate leads to $H^2$-bounds for $\chi(u_h(\lambda)-u_h(0))$ on either side of the boundary $\partial X_\epsilon$, locally uniform in $\lambda$.  Taking the jump of the normal derivatives is continuous into $H^{1/2}(\partial X_\epsilon)$, so
$\lambda\mapsto R_h(\lambda)-R_h(0)$
is analytic from $H^1(\partial X_\epsilon)$ to $H^{1/2}(\partial X_\epsilon)$, and hence to $L^2(\partial X_\epsilon)$. Adding the fixed bounded operator $R_h(0)$ proves the analyticity of $R_h(\lambda)$  asserted in~\eqref{eq3.4}.

Now, let us prove the analyticity of the difference $R_{\widetilde m}(\lambda)-R_m(\lambda)$ asserted in~\eqref{eq3.5}. 
In what follows, $\Delta_{\mathcal U,\widetilde m}=\bigoplus_j\Delta_{\mathcal U^j,\widetilde m}$ stands for  the Dirichlet Laplacian on $\mathcal U$. 
  By
$\mathcal P_{\mathcal U,h}(\lambda)$ we denote the Poisson operator on
$\mathcal U$ for the metric $h\in\{m,\widetilde m\}$, defined by
$$
 \mathcal P_{\mathcal U,h}(\lambda)\varphi
 =u_h(\lambda)|_{\mathcal U},
$$
where $u_h(\lambda)$ is the weak solution of the Dirichlet  problem ~\eqref{eq3.1}.

For $\varphi\in C^\infty(\partial X_\epsilon)$, the function
$\chi\mathcal P_{\mathcal U,m}(\lambda)\varphi$ has trace
$\varphi$ on $\partial X_\epsilon$, vanishes near the conical points, and in the interior of
each disk satisfies
\begin{equation}\label{eq3.11}
 (\Delta_{\widetilde m}+\lambda)
 \bigl(\chi\mathcal P_{\mathcal U,m}(\lambda)\varphi\bigr)
 =(\Delta_m+\lambda)
 \bigl(\chi\mathcal P_{\mathcal U,m}(\lambda)\varphi\bigr)=
 [\Delta_m,\chi]\mathcal P_{\mathcal U,m}(\lambda)\varphi.
\end{equation}
Note that the difference
$\mathcal P_{\mathcal U,\widetilde m}(\lambda)\varphi
-\chi\mathcal P_{\mathcal U,m}(\lambda)\varphi$ has zero boundary trace.
Applying the Dirichlet resolvent to the difference and taking into account~\eqref{eq3.11} give
\begin{equation}\label{eq3.12}
 \mathcal P_{\mathcal U,\widetilde m}(\lambda)
 -
 \chi\mathcal P_{\mathcal U,m}(\lambda)=
 -
 (\Delta_{\mathcal U,\widetilde m}+\lambda)^{-1}
 [\Delta_m,\chi]\mathcal P_{\mathcal U,m}(\lambda).
\end{equation}
Since $\chi=1$ near $\partial X_\epsilon$,  by  taking the
normal derivative  in \eqref{eq3.12} at $\partial X_\epsilon$, we obtain
\begin{equation}\label{eq3.13}
 R_{\widetilde m}(\lambda)-R_m(\lambda)
 =
 -\partial_\nu
 (\Delta_{\mathcal U,\widetilde m}+\lambda)^{-1}
 [\Delta_m,\chi]\mathcal P_{\mathcal U,m}(\lambda);
\end{equation}
it is important here that  for both  metrics the Dirichlet problem on $X_\epsilon$ remains the same.
We rewrite
\eqref{eq3.13} as follows:
\begin{equation}\label{eq3.14}
 R_{\widetilde m}(\lambda)-R_m(\lambda)
 =
 -\Bigl(
 \partial_\nu
 (\Delta_{\mathcal U,\widetilde m}+\lambda)^{-1}\eta
 \Bigr)
 [\Delta_m,\chi]
 \Bigl(
 \eta\mathcal P_{\mathcal U,m}(\lambda)
 \Bigr).
\end{equation}
Recall that  $\eta$ stands for the cutoff function introduced in the beginning of the proof.
The support of $\eta$ is separated from both
$\partial X_\epsilon$ and the conical points. The standard local smooth elliptic regularity theory therefore implies that the left and the right factors 
$$
S_\delta\ni\lambda\longmapsto \Bigl(
 \partial_\nu
 (\Delta_{\mathcal U,\widetilde m}+\lambda)^{-1}\eta
 \Bigr)\in \mathcal B\bigl(
 L^2(\mathcal U),
 H^2(\partial X_\epsilon)
 \bigr),
$$
$$
S_\delta\ni\lambda\longmapsto \eta\mathcal P_{\mathcal U,m}(\lambda)\in \mathcal B\bigl(
 H^1(\partial X_\epsilon),
 H^1(\mathcal U)
 \bigr)
$$
in~\eqref{eq3.14} are analytic families of bounded operators in the listed spaces.
Since the commutator  $[\Delta_m,\chi]$   is a $\lambda$-independent continuous operator from $H^1(\mathcal U)$ to $L^2(\mathcal U)$,
 the factorization~\eqref{eq3.14} proves analyticity of the family $ R_{\widetilde m}(\lambda)-R_m(\lambda)$ in~\eqref{eq3.5}.

Since the Dirichlet Laplacians involved are strictly positive, their resolvents are analytic in a neighborhood of zero. The preceding arguments therefore establish analyticity of the families~\eqref{eq3.4} and~\eqref{eq3.5}  in a neighborhood of zero as well. This completes the proof of the first assertion. 
%======

We now prove the second assertion. 
Let $\lambda\in S_\delta$ and $h\in\{m,\widetilde m\}$.
If $R_h(\lambda)\varphi=0$, the solution $u_h(\lambda)$ of the problem~\eqref{eq3.1}
extends across the cut to a global weak solution of
$(\Delta_h+\lambda)u=0$. Nonnegativity of the Laplacian $\Delta_h$
and $\Re\lambda>0$ imply that $u=0$, and hence $\varphi=0$.
The spectrum of $R_h(\lambda)$ consists of isolated eigenvalues of finite
algebraic multiplicity. Since its kernel is trivial,
zero belongs to its resolvent set.
As is well-known~\cite{Ref20}, this together with analyticity of  $\lambda\mapsto R_h(\lambda)$ in~\eqref{eq3.4} implies analyticity of
the operator family
$$
 S_\delta\ni\lambda\longmapsto R_h(\lambda)^{-1}
 \in\mathcal B\bigl(
 L^2(\partial X_\epsilon),
 H^1(\partial X_\epsilon)
 \bigr).
$$
Hence, due to the representation~\eqref{eq3.2} for $T(\lambda)-I$ and  analyticity of $R_{\widetilde m}(\lambda)-R_m(\lambda)$ in~\eqref{eq3.5},  the family
$$
 S_\delta\ni\lambda\longmapsto T(\lambda)-I
 \in\mathcal B\bigl(
 L^2(\partial X_\epsilon),
 H^2(\partial X_\epsilon)
 \bigr)
$$
is also analytic. As the boundary $\partial X_\epsilon$ is one-dimensional, the embedding
of 
 $H^2(\partial X_\epsilon)$
into
 $L^2(\partial X_\epsilon)$
is trace class; e.g.~\cite{Ref27}. This leads to  the analyticity of $\lambda\mapsto T(\lambda)-I$ in the trace norm as asserted
in~\eqref{eq3.6}.
%======

Before deducing the decay estimate~\eqref{eq3.7} for $ T(\lambda)-I$,  we  derive a suitable representation for this operator.

Let
$ u_m(\lambda)=(\Delta_m+\lambda)^{-1}
(v\otimes\delta_{\partial X_\epsilon})$, where  $v\in L^2(\partial X_\epsilon)$. 
Here the distribution $v\otimes\delta_{\partial X_\epsilon}$  acts on a function $\phi$ in the 
form domain $ \mathcal V=H^1(X, m)$ of $\Delta_m$ by
$\int_{\partial X_\epsilon}v{\left.\overline\phi\right|_{\partial X_\epsilon}}\,ds$.
By the trace theorem we have
\begin{equation}\label{eq3.15}
\left|\int_{\partial X_\epsilon}
v{\left.\overline\phi\right|_{\partial X_\epsilon}}\,ds\right|
\leq C\|v\|_{L^2(\partial X_\epsilon)}
\bigl(\|d\phi\|_{L^2(X,m)}^2
+\|\phi\|_{L^2(X,m)}^2\bigr)^{1/2},
\end{equation}
where the last factor is the norm $\|\phi\|_{\mathcal V}$ of $\phi$ in $\mathcal V$.  Thus $v\otimes\delta_{\partial X_\epsilon}$
belongs to the dual space $\mathcal V^*=H^{-1}(X, m)$, and
$u_m(\lambda)$
is the unique weak solution in $\mathcal V$ to the problem 
\begin{equation*}
(\Delta_m+\lambda)u_m(\lambda)=v\otimes\delta_{\partial X_\epsilon}\quad\text{ on }\quad X.
\end{equation*}
As is well-known, the inverse of the Neumann jump can be written in the  form
$$
 R_m(\lambda)^{-1}v
 =\left.(\Delta_m+\lambda)^{-1}
 (v\otimes\delta_{\partial X_\epsilon})
 \right|_{\partial X_\epsilon},
$$
e.g.~\cite{Ref6}. 
Thus the function $u_m(\lambda)$ is also the weak solution to the problem~\eqref{eq3.1} with $h=m$ and  $  \varphi=  R_m(\lambda)^{-1}v$. 
 Unique solvability of the
Dirichlet problems in the disks then implies
\begin{equation}\label{eq3.16}
 \mathcal P_{\mathcal U,m}(\lambda)
 R_m(\lambda)^{-1}v
 =  \left. u_m(\lambda)\right|_{\mathcal U}=
 \left.(\Delta_m+\lambda)^{-1}
 (v\otimes\delta_{\partial X_\epsilon})
 \right|_{\mathcal U}.
\end{equation}
Now, thanks to the equalities~\eqref{eq3.13} and~\eqref{eq3.16} together with the representation \eqref{eq3.2} for the relative jump $T(\lambda)$, we arrive at the representation
$$
 T(\lambda)-I= -\partial_\nu
 (\Delta_{\mathcal U,\widetilde m}+\lambda)^{-1}
 [\Delta_m,\chi](\Delta_m+\lambda)^{-1}
 (\,\cdot\,\otimes\delta_{\partial X_\epsilon}).
$$
It is convenient to rewrite this representation in the form
\begin{equation}\label{eq3.17}
 T(\lambda)-I
 =-\left(
 \partial_\nu
 (\Delta_{\mathcal U,\widetilde m}+\lambda)^{-1}\eta
 \right)\,[\Delta_m,\chi]\left(
 \eta(\Delta_m+\lambda)^{-1}
 (\,\cdot\,\otimes\delta_{\partial X_\epsilon})
 \right),
\end{equation}
where  $\eta$ is the cutoff function introduced in the beginning of the proof. For ease of reference in the estimates below, we distinguish three factors on the right hand side of~\eqref{eq3.17}: the left resolvent factor, the commutator, and the right resolvent factor.

We next establish the estimate~\eqref{eq3.7}.
The trace estimate~\eqref{eq3.15} shows that the mapping  $L^2(\partial X_\epsilon)\ni v \mapsto v\otimes\delta_{\partial X_\epsilon}\in
\mathcal V^*$ is continuous.  Since
$$\|\cdot\|_{\mathcal V}
=\|(\Delta_m+1)^{1/2}\cdot\|_{L^2(X,m)},\quad 
\|\cdot\|_{\mathcal V^*}
=\|(\Delta_m+1)^{-1/2}\cdot\|_{L^2(X,m)},
$$
the spectral theorem implies
$$
 \| (\Delta_m+\lambda)^{-1}\|_{\mathcal B(\mathcal V^*,\mathcal V)}
 =
 \sup_{\mu\in\sigma(\Delta_m)}
 \frac{\mu+1}{|\mu+\lambda|}
 \leq C,
 \qquad \lambda\in S_\delta,\quad |\lambda|\geq1,
$$
where $C$ depends on $\delta$.
Multiplication by the cutoff function $\eta$ is bounded from $\mathcal V$ to
$H^1(\mathcal U)$.  Consequently, for the right resolvent factor in~\eqref{eq3.17} we obtain
\begin{equation}\label{eq3.18}
 \left\|
 \eta(\Delta_m+\lambda)^{-1}  (\,\cdot\,\otimes\delta_{\partial X_\epsilon})
 \right\|_{\mathcal B(L^2(\partial X_\epsilon),H^1(\mathcal U))}
 \leq C,
 \qquad \lambda\in S_\delta,\quad |\lambda|\geq1.
\end{equation}

%======-=-================

To estimate the left resolvent factor in~\eqref{eq3.17}, 
we shall rely on the standard  elliptic  parameter-dependent  estimate
\begin{equation}\label{eq3.19}
 \begin{aligned}
 \|u\|_{H^2(\mathcal U)}
 +|\lambda|^{1/2}\|u\|_{H^1(\mathcal U)}
 +|\lambda|\|u\|_{L^2(\mathcal U)}
 \leq C\|(\Delta_{\widetilde m}+\lambda)u\|_{L^2(\mathcal U)},
 \end{aligned}
\end{equation}
where  $C$ is independent of
$u\in H^2(\mathcal U)\cap H^1_0(\mathcal U)$ and
$\lambda\in S_\delta$ with $|\lambda|\geq\lambda_0\geq 1$;
see e.g.~\cite{Ref1} or~\cite[Sec.~3.6.3]{Ref24}.

Let $f\in L^2(\mathcal U)$ and $\lambda\in S_\delta$. Then for 
$$
 w_\lambda=(\Delta_{\mathcal U,\widetilde m}+\lambda)^{-1}\eta f
$$
 the estimate~\eqref{eq3.19} with $u=w_\lambda$ gives \begin{equation}\label{eq3.20}
 \|w_\lambda\|_{H^1(\mathcal U)}
 \leq C|\lambda|^{-1/2}\|f\|_{L^2(\mathcal U)},
 \qquad\lambda\in S_\delta,\quad|\lambda|\geq\lambda_0.
\end{equation}

Within a   neighborhood of the boundary $\partial\mathcal U=\partial X_\epsilon$ in
$\mathcal U$, choose a chain of neighborhoods 
$$
 V_0\subset V_1\subset\cdots\subset V_{2n+4},
 \qquad \overline{V_i}\subset V_{i+1},
 \quad i=0,\ldots,2n+3,\quad \overline{V_{2n+4}}\cap\operatorname{supp}\eta=\varnothing,
 $$
where the closures are taken in $\mathcal U$.  Choose a smooth cutoff function $\rho_i$ on $\mathcal U$, supported in $V_{i+1}$ and equal
to one near $\overline{V_i}$.  The function $\rho_iw_\lambda$ in the domain of the Dirichlet Laplacian $\Delta_{\mathcal U,\widetilde m}$ satisfies
$$
 (\Delta_{\mathcal U,\widetilde m}+\lambda)(\rho_iw_\lambda)
 =[\Delta_{\mathcal U,\widetilde m},\rho_i]w_\lambda.
$$
Since the commutator is a first order differential operator with smooth coefficients, substitution of $u=\rho_iw_\lambda$ into 
\eqref{eq3.19} implies 
\begin{equation}\label{eq3.21}
\|w_\lambda\|_{H^2(V_i)}
 +|\lambda|^{1/2}\|w_\lambda\|_{H^1(V_i)}
 +|\lambda|\|w_\lambda\|_{L^2(V_i)}\leq C_i\|w_\lambda\|_{H^1(V_{i+1})},
 \end{equation}
and, as a consequence,
\begin{equation}\label{eq3.22}
 \|w_\lambda\|_{H^1(V_i)}
 \leq C_i |\lambda|^{-1/2}
 \|w_\lambda\|_{H^1(V_{i+1})},
 \qquad\lambda\in S_\delta, \quad |\lambda|\geq\lambda_0.
\end{equation}

Starting from~\eqref{eq3.20}, we iterate
\eqref{eq3.22} over the $2n+2$ pairs from
$V_{2n+4}$ to $V_2$, gaining the factor $|\lambda|^{-n-1}$.
The $H^2$-term of~\eqref{eq3.21}
with $i=1$ then yields
\begin{equation}\label{eq3.23}
 \|w_\lambda\|_{H^2(V_1)}
 \leq C_n |\lambda|^{-n-3/2}\|f\|_{L^2(\mathcal U)}.
\end{equation}
Since $\Delta_{\widetilde m}w_\lambda=-\lambda w_\lambda$ in $V_1$,
the standard local elliptic a priori estimate
and~\eqref{eq3.23}, with
$|\lambda|\geq\lambda_0\geq1$, yield
\begin{equation}\label{eq3.24}
 \begin{aligned}
 \|w_\lambda\|_{H^4(V_0)}
 &\leq C\bigl(
 \|\Delta_{\widetilde m}w_\lambda\|_{H^2(V_1)}
 +\|w_\lambda\|_{L^2(V_1)}
 \bigr)\\
 &\leq C|\lambda|\|w_\lambda\|_{H^2(V_1)}
 \leq C_n|\lambda|^{-n-1/2}\|f\|_{L^2(\mathcal U)}.
 \end{aligned}
\end{equation}
By the boundary trace theorem and the continuous embedding
$H^{5/2}\hookrightarrow H^2$ on $\partial X_\epsilon$,
the estimate~\eqref{eq3.24} implies
\begin{equation*}
 \|\partial_\nu w_\lambda\|_{H^2(\partial X_\epsilon)}
 \leq C\|w_\lambda\|_{H^4(V_0)}\\
 \leq C_n|\lambda|^{-n-1/2}\|f\|_{L^2(\mathcal U)},
 \quad\lambda\in S_\delta, \quad |\lambda|\geq\lambda_0.
 \end{equation*}
By this estimate and the trace-class embedding
$H^2(\partial X_\epsilon)
\hookrightarrow L^2(\partial X_\epsilon)$,
the left resolvent factor satisfies
\begin{equation}\label{eq3.25}
 \bigl\|\partial_\nu
 (\Delta_{\mathcal U,\widetilde m}+\lambda)^{-1}\eta\bigr\|_{\mathfrak S_1\bigl(L^2(\mathcal U),\,L^2(\partial X_\epsilon)\bigr)}
 \leq C_n |\lambda|^{-n-1/2},
 \qquad\lambda\in S_\delta, \quad |\lambda|\geq\lambda_0.
\end{equation}

In the factorization~\eqref{eq3.17}, the commutator
is bounded from $H^1(\mathcal U)$ to $L^2(\mathcal U)$.
The estimates~\eqref{eq3.18}
and~\eqref{eq3.25} therefore imply
\eqref{eq3.7} by the two-sided ideal property
of the trace class in the algebra of bounded operators.

Finally,  we prove the third assertion. The commutation relations between $R_h(\lambda)$ and the
orthogonal projectors $\Pi^\sharp$ extend from $\lambda>0$ to
$S_\delta$ by analyticity, so the representations~\eqref{eq3.3} hold
throughout $S_\delta$. Since the projectors are independent of
$\lambda$, the restrictions $T^\sharp(\lambda)-I$ are analytic
in the trace norm. The estimate~\eqref{eq3.7} then yields
$$
 \|T^\sharp(\lambda)-I\|_{\mathfrak S_1}
 \leq\|T(\lambda)-I\|_{\mathfrak S_1}
 \leq C_n|\lambda|^{-n},
 \qquad
 \lambda\in S_\delta,\quad|\lambda|\geq\lambda_0,
 \quad\sharp\in\{D,N\}.
$$
Continuity of the Fredholm determinant
implies~\eqref{eq3.8}.

The operators $R_m(\lambda)$ and $R_{\widetilde m}(\lambda)$
are invertible and commute with $\Pi^\sharp$.
Their restrictions $R_h^\sharp$ are therefore invertible, and so is
$T^\sharp(\lambda)$. Hence its Fredholm determinant is
analytic and nonzero on $S_\delta$.

By the standard continuity estimate for the Fredholm determinant
\cite[(3.7)]{Ref28} and the decay estimate~\eqref{eq3.7} for $T(\lambda)-I$,
for sufficiently large $|\lambda|$  we have
$$
 \begin{aligned}
 \bigl|\log\det\nolimits_F T^\sharp(\lambda)\bigr|
 &\leq C\bigl|\det\nolimits_F T^\sharp(\lambda)-1\bigr|\\
 &\leq C\|T^\sharp(\lambda)-I\|_{\mathfrak S_1} e^{1+\|T^\sharp(\lambda)-I\|_{\mathfrak S_1}}
 \leq C_n|\lambda|^{-n},\quad \lambda\in S_\delta.
 \end{aligned}
$$
This establishes the decay estimate~\eqref{eq3.9} for $\log\det_F T^\sharp(\lambda)$ and completes the proof.
\end{proof}
%%=====
%%=====

\section{Relative BFK identity}\label{Sec4}

For $h\in\{m,\widetilde m\}$, let $\Delta_h$ be the Friedrichs
Laplacian on $(X,h)$, and let
$\Delta_{\partial X_\epsilon,h}$ be the Friedrichs
Dirichlet Laplacian on
$X\setminus\partial X_\epsilon$.
The pullback $\iota^*$ by $\iota$ is unitary in the corresponding
$L^2$-spaces; then  $ \Pi^D=\frac12(I-\iota^*)$ and 
$ \Pi^N=\frac12(I+\iota^*)$
are orthogonal projections in these spaces.  Both Laplacians
commute with $\Pi^D$ and $\Pi^N$.  Denote the restrictions of
$\Delta_h$ to
$\Pi^D L^2(X,h)$ and
$\Pi^N L^2(X,h)$ by $\Delta_h^D$ and
$\Delta_h^N$, respectively.  On either half of the double these restrictions are
unitarily equivalent to the natural Dirichlet and Neumann Laplacians
defined in the introduction.  For
$\sharp\in\{D,N\}$, let $\Delta_{\partial X_\epsilon,h}^\sharp$ denote the
corresponding restriction of
$\Delta_{\partial X_\epsilon,h}$.

For a single metric $h$, the BFK gluing formula~\cite{Ref5, Ref11} has the
form
$$
 \log\Det(\Delta_h+\lambda)
 -\log\Det(\Delta_{\partial X_\epsilon,h}+\lambda)
 =
 \log\Det R_h(\lambda)+C_h,
$$
where $C_h$ is independent of $\lambda$. Its traditional proof first establishes this formula up to a constant
$C_h$.  The value of $C_h$ is then found from the separate
large-$\lambda$ asymptotics of the two Laplacian determinants and
the Neumann jump determinant. Applying this argument after
restriction to $\Pi^D$ and $\Pi^N$ would require additional
equivariant asymptotics.

Instead, in the proof of Proposition~\ref{Prop4.1},   we introduce a relative zeta function and express it as a contour integral
of the differences of the Laplacian resolvents.  In the difference of the zeta functions for $\Delta_{\partial X_\epsilon,m}^\sharp$ and $\Delta_{\partial X_\epsilon,\widetilde m}^\sharp$,
the common summand, corresponding to the Dirichlet Laplacian on $X_\epsilon$ with $m=\widetilde m$, cancels before analytic continuation. This
leaves only the zeta functions of the Laplacians on the disks studied in Section~\ref{Sec2}. In particular, no study of  individual  zeta functions of
 $\Delta_{\partial X_\epsilon,m}^\sharp$ or $\Delta_{\partial X_\epsilon,\widetilde m}^\sharp$ at zero is required.
Lemma~\ref{Lem3.1} then provides the
required relative normalization at $\lambda$-infinity and fixes the remaining
constant.

\begin{proposition}[Relative BFK identity]\label{Prop4.1}
For every $\lambda>0$ and $\sharp\in\{D,N\}$, we have
\begin{equation}\label{eq4.1}
 \frac{\Det(\Delta_{\widetilde m}^\sharp+\lambda)}
      {\Det(\Delta_m^\sharp+\lambda)}
 \frac{\Det(\Delta_{\mathcal U,m}^\sharp+\lambda)}
      {\Det(\Delta_{\mathcal U,\widetilde m}^\sharp+\lambda)}
 =
 \det\nolimits_F T^\sharp(\lambda),
\end{equation}
where $
\Delta_{\mathcal U,h}^\sharp
=\oplus_j\Delta_{\mathcal U^j,h}^\sharp$,  the determinants of the Laplacians are the spectral zeta regularized determinants, and $\det\nolimits_F T^\sharp(\lambda)$ is the Fredholm determinant of the relative Neumann jump operator $T^\sharp(\lambda)$ from \eqref{eq3.3} and Lemma~\ref{Lem3.1}.
\end{proposition}

\begin{proof}

We first derive the resolvent comparison identity
\begin{equation}\label{eq4.2}
 (\Delta_h^\sharp+\lambda)^{-1}
 -(\Delta_{\partial X_\epsilon,h}^\sharp+\lambda)^{-1}
 =
 \mathcal P_h(\lambda)
 \bigl(R_h^\sharp(\lambda)\bigr)^{-1}
 \bigl(\mathcal P_h(\lambda)\bigr)^*,\quad h\in\{m,\widetilde m\}.
\end{equation}
Here  
$\mathcal P_h(\lambda)$ is the Poisson operator on $ X\setminus  \partial X_\epsilon$, defined by
$ \mathcal P_h(\lambda)\varphi=u_h(\lambda)$,
where $u_h(\lambda)$ is the weak solution to the Dirichlet problem~\eqref{eq3.1}.

With this aim in mind we note that  the difference
$$
u:= (\Delta_h^\sharp+\lambda)^{-1}f
 -(\Delta_{\partial X_\epsilon,h}^\sharp+\lambda)^{-1}f, \qquad f\in \Pi^\sharp L^2(X, h),$$
satisfies the boundary value problem
$$
(\Delta_h+\lambda)u=0\quad \text{in}\quad X\setminus  \partial X_\epsilon;  \qquad   u=(\Delta_h^\sharp+\lambda)^{-1}f  \quad \text{on} \quad \partial X_\epsilon.
$$
As a consequence, we have
 \begin{equation}\label{eq4.3}
 (\Delta_h^\sharp+\lambda)^{-1}f
 -(\Delta_{\partial X_\epsilon,h}^\sharp+\lambda)^{-1}f= \mathcal P_h (\lambda)\left.\left((\Delta_h^\sharp+\lambda)^{-1}f\right)\right|_{\partial X_\epsilon}.
 \end{equation}
The argument leading to
\eqref{eq3.16}, now applied to $h$ and  $X\setminus  \partial X_\epsilon$, 
 gives
$$
 (\Delta_h^\sharp+\lambda)^{-1}(v\otimes\delta_{\partial X_\epsilon})
 =
 \mathcal P_h(\lambda)
 \bigl(R_h^\sharp(\lambda)\bigr)^{-1}v.
$$
Both sides are bounded operators from $\Pi^\sharp L^2(\partial X_\epsilon)$ to $\Pi^\sharp L^2(X,h)$.
Passing to the adjoint operators yields
$$
\left.\left( (\Delta_h^\sharp+\lambda)^{-1}f \right)\right|_{\partial X_\epsilon}
 =
 \bigl(R_h^\sharp(\lambda)\bigr)^{-1}
 \mathcal P_h^*(\lambda)f.
$$
This together with~\eqref{eq4.3} implies the resolvent comparison identity~\eqref{eq4.2}.

Next, we show that the left hand side of the comparison identity~\eqref{eq4.2} is trace class in $\Pi^\sharp L^2(X,h)$ and express the difference between its traces for $h=m$ and $h=\widetilde m$ in terms of the Fredholm determinant of the relative Neumann jump. 

Let
$$
 u=\mathcal P_h(\lambda)\varphi,
 \qquad
 v=\mathcal P_h(\bar\lambda)\psi,
 \qquad
 \varphi,\psi\in\Pi^\sharp C^\infty(\partial X_\epsilon).
$$
Differentiation of the boundary value problem for $u$ gives
$$
 (\Delta_h+\lambda)\partial_\lambda u=-u,
 \qquad
 \left.\partial_\lambda u\right|_{\partial X_\epsilon}=0.
$$
Hence the Green formula for this problem takes the form
$$
 -\langle u,v\rangle_{L^2(X,h)}
 =
 -\left\langle
  \partial_\lambda R_h^\sharp\varphi,\psi
 \right\rangle_{L^2(\partial X_\epsilon)}.
$$
Therefore,
$$
 \left\langle
  \partial_\lambda R_h^\sharp(\lambda)\varphi,\psi
 \right\rangle_{L^2(\partial X_\epsilon)}
 =
 \left\langle
  \mathcal P_h(\lambda)\varphi,
  \mathcal P_h(\bar\lambda)\psi
 \right\rangle_{L^2(X,h)},
$$
or, equivalently, 
\begin{equation}\label{eq4.4}
 \partial_\lambda R_h^\sharp(\lambda)
 =
 \left.
 \mathcal P_h^*(\bar\lambda)\mathcal P_h(\lambda)
 \right|_{\Pi^\sharp L^2(\partial X_\epsilon)}.
\end{equation}

Recall that $R_h(\lambda)$, $\lambda>0$, is positive and selfadjoint, and notice that  the operator
$$
 (R_h(\lambda))^{-1/2}(\partial_\lambda R_h(\lambda))(R_h(\lambda))^{-1/2}
$$
is trace class in
$L^2(\partial X_\epsilon)$ as a classical pseudodifferential operator of order $-2$;  see, e.g.,~\cite{Ref27}.
Since $R_h$ and $\partial_\lambda R_h$ commute with
the orthogonal projections $\Pi^\sharp$, we have
$$
 \bigl(R_h^\sharp(\lambda)\bigr)^{-1/2}
 \bigl(\partial_\lambda R_h^\sharp(\lambda)\bigr)
 \bigl(R_h^\sharp(\lambda)\bigr)^{-1/2}
 =\left.
 (R_h)^{-1/2}(\partial_\lambda R_h)(R_h)^{-1/2}
 \right|_{\Pi^\sharp L^2(\partial X_\epsilon)}.
$$
The operator on the left is therefore trace class too.  By
\eqref{eq4.4}, it equals
$$
 (\mathcal P_h(R_h^\sharp)^{-1/2})^*
 (\mathcal P_h(R_h^\sharp)^{-1/2}).
$$
Consequently, $\mathcal P_h(R_h^\sharp)^{-1/2}$ is
Hilbert--Schmidt. This implies that  the right hand side of the resolvent identity~\eqref{eq4.2} is trace class, and the 
trace is
\begin{equation}\label{eq4.5}
\begin{aligned}
 \Tr_{\Pi^\sharp L^2(X,h)}\!\left[
  (\Delta_h^\sharp+\lambda)^{-1}
  -(\Delta_{\partial X_\epsilon,h}^\sharp+\lambda)^{-1}
 \right]
 \\=   \Tr_{\Pi^\sharp L^2(X,h)}\!\left[ (\mathcal P_h(R_h^\sharp)^{-1/2}) (\mathcal P_h(R_h^\sharp)^{-1/2})^*\right]
 \\= \Tr_{\Pi^\sharp L^2(\partial X_\epsilon)}\!\left[ (\mathcal P_h(R_h^\sharp)^{-1/2})^*
 (\mathcal P_h(R_h^\sharp)^{-1/2}) 
 \right]
 \\=
 \Tr_{\Pi^\sharp L^2(\partial X_\epsilon)}\!\left[
  \bigl(\partial_\lambda R_h^\sharp(\lambda)\bigr)
  \bigl(R_h^\sharp(\lambda)\bigr)^{-1}
 \right],
 \end{aligned}
\end{equation}
where the last equality follows from
\eqref{eq4.4} and the cyclicity of the trace.
Let us note that the trace identity~\eqref{eq4.5}, without restriction to functions that
are even or odd with respect to the reflection $\iota$, is well
known, e.g.~\cite[Theorem~2.2]{Ref6}.
 
By definition of the relative Neumann jump operator we have $
 R_{\widetilde m}^\sharp
 =
 T^\sharp R_m^\sharp
$, cf.~\eqref{eq3.3}.
Differentiating this identity and multiplying on the right by
$
 \bigl(R_{\widetilde m}^\sharp\bigr)^{-1}$
we obtain
$$
 \bigl(\partial_\lambda R_{\widetilde m}^\sharp\bigr)
 \bigl(R_{\widetilde m}^\sharp\bigr)^{-1}
 =
 \bigl(\partial_\lambda T^\sharp\bigr)
 \bigl(T^\sharp\bigr)^{-1}
 +
 T^\sharp
 \bigl(\partial_\lambda R_m^\sharp\bigr)
 \bigl(R_m^\sharp\bigr)^{-1}
 \bigl(T^\sharp\bigr)^{-1}.
$$
Here the last term has the same trace as
$(\partial_\lambda R_m^\sharp)(R_m^\sharp)^{-1}$, by cyclicity
of the trace for products of a trace-class operator and bounded
operators.

 Subtracting the two traces and expressing each of them in terms of resolvents as in~\eqref{eq4.5}, we arrive at
\begin{equation}\label{eq4.6}
 \begin{aligned}
\Tr_{\Pi^\sharp L^2(X,\widetilde m)}\Bigl[
  (&\Delta_{\widetilde m}^\sharp+\lambda)^{-1}
  -(\Delta_{\partial X_\epsilon,\widetilde m}^\sharp+\lambda)^{-1}
 \Bigr]\\
 &-\Tr_{\Pi^\sharp L^2(X,m)}\!\left[
  (\Delta_m^\sharp+\lambda)^{-1}
  -(\Delta_{\partial X_\epsilon,m}^\sharp+\lambda)^{-1}
 \right]
 \\
 &=
 \Tr_{\Pi^\sharp L^2(\partial X_\epsilon)}\!\left[(\partial_\lambda T^\sharp(\lambda))
             T^\sharp(\lambda)^{-1}\right]
 =
 \partial_\lambda\log\det\nolimits_F T^\sharp(\lambda).
 \end{aligned}
 \end{equation}
 Thus we expressed the difference between the traces of the left hand side of~\eqref{eq4.2}  for $h=m$ and $h=\widetilde m$ in terms of the Fredholm determinant of the relative Neumann jump operator $T^\sharp$. 
%=========

For the conical metric $h=m$, the short-time heat trace estimates
$
 \Tr e^{-t\Delta_m}=O(t^{-1})$  and 
$\Tr e^{-t\Delta_{\partial X_\epsilon,m}}
 =O(t^{-1})$   as $t\to0^+$,
imply absolute convergence of the spectral series for
$(\Delta_m+\lambda)^{-s}$ and
$(\Delta_{\partial X_\epsilon,m}+\lambda)^{-s}$
when $\lambda>0$ and $\Re s>1$;
see~\cite[Lemma~2.5 and (2.13)]{Ref11}. For the smoothed metric
$h=\widetilde m$, the same conclusion follows from the standard
smooth theory. Since the Laplacians commute with the orthogonal
projectors $\Pi^D$ and $\Pi^N$, the odd and even spectral
subseries also converge absolutely in this half-plane.
We may therefore define the spectral zeta functions
$$
 \zeta_h^\sharp(s;\lambda)
 :=
 \Tr(\Delta_h^\sharp+\lambda)^{-s},
 \qquad
 \zeta_{\partial X_\epsilon,h}^\sharp(s;\lambda)
 :=
 \Tr(\Delta_{\partial X_\epsilon,h}^\sharp
      +\lambda)^{-s}, \quad \lambda>0,\quad \Re s>1.
$$

%=======

We now obtain a contour representation of the relative zeta function
\begin{equation}\label{eq4.7}
 \zeta_{\mathrm{rel}}^\sharp(s;\lambda)
 :=
\left( \zeta_{\widetilde m}^\sharp(s;\lambda) -\zeta_{\partial X_\epsilon,\widetilde m}^\sharp(s;\lambda)\right)
-\left(\zeta_m^\sharp(s;\lambda)
 -\zeta_{\partial X_\epsilon,m}^\sharp(s;\lambda)\right), \quad  \Re s>1,
\end{equation}
and continue it analytically to $s=0$.

 Let $\mathcal C_\lambda$ be a contour in $S_\delta$ running clockwise around the cut $[\lambda,\infty)$, where $\lambda>0$. Outside a bounded set, the contour consists of two horizontal half-lines.
By the contour formula for complex powers we have
\begin{equation}\label{eq4.8}
\begin{aligned}
 & (\Delta_h^\sharp+\lambda)^{-s}
 -(\Delta_{\partial X_\epsilon,h}^\sharp+\lambda)^{-s}\\
 &\quad=
 \frac1{2\pi i}\int_{\mathcal C_\lambda}
 (\lambda-\mu)^{-s}
 \bigl((\Delta_h^\sharp+\mu)^{-1}
 -(\Delta_{\partial X_\epsilon,h}^\sharp+\mu)^{-1}\bigr)
 \,d\mu,\quad \Re s>1,
\end{aligned}
\end{equation}
where we use the principal branch of
$(\lambda-\mu)^{-s}$, with cut $[\lambda,\infty)$; cf.~\cite{Ref27}.
Next, we show that this integral converges in the trace norm, which permits
taking the trace under the integral.

Indeed, the resolvent identity takes the form
$$
\begin{aligned}
 & (\Delta_h^\sharp+\mu)^{-1}
 -(\Delta_{\partial X_\epsilon,h}^\sharp+\mu)^{-1}\\
 &=
 \bigl(I+(1-\mu)(\Delta_h^\sharp+\mu)^{-1}\bigr)
 \bigl((\Delta_h^\sharp+1)^{-1}
 -(\Delta_{\partial X_\epsilon,h}^\sharp+1)^{-1}\bigr)
 \bigl(I+(1-\mu)
 (\Delta_{\partial X_\epsilon,h}^\sharp+\mu)^{-1}\bigr).
\end{aligned}
$$
The factor in the middle is trace class as before. For $\mu\in\mathcal C_\lambda$ both bordering factors are uniformly bounded in the operator norm as the distance from $-\mu$ to the spectrum of either Laplacian is at least $|\mu|$ and $\inf_{\mu\in\mathcal C_\lambda}|\mu|>0$.
 The two-sided ideal
property of the trace class therefore implies
$$
 \bigl\|(\Delta_h^\sharp+\mu)^{-1}
 -(\Delta_{\partial X_\epsilon,h}^\sharp+\mu)^{-1}
 \bigr\|_{\mathfrak S_1}
 \le
 C(\lambda)
 \bigl\|(\Delta_h^\sharp+1)^{-1}
 -(\Delta_{\partial X_\epsilon,h}^\sharp+1)^{-1}
 \bigr\|_{\mathfrak S_1},
 \quad \mu\in\mathcal C_\lambda,
$$
where the trace norms are taken in
$\Pi^\sharp L^2(X,h)$ and $C(\lambda)$ does not depend on $\mu$.

Clearly, the  contour integral
$
\int_{\mathcal C_\lambda}|(\lambda-\mu)^{-s}|\,|d\mu|
$
converges for $\Re s>1$. The preceding estimate therefore implies convergence of the integral~\eqref{eq4.8} in the trace norm. Continuity of the trace in that norm justifies taking the trace under the integral.

On the next step, we rely on the trace identity~\eqref{eq4.6} extended to the sector $S_\delta$ by analyticity of both sides.
Namely, subtracting the resulting trace formulas for $h=\widetilde m$
and $h=m$, we obtain
$$
 \zeta_{\mathrm{rel}}^\sharp(s;\lambda)
 =
 \frac1{2\pi i}\int_{\mathcal C_\lambda}
 (\lambda-\mu)^{-s}
 \,\partial_\mu\log\det\nolimits_F T^\sharp(\mu)\,d\mu,
 \qquad \Re s>1.
$$

To compute the derivative at $s=0$, we integrate by parts. The boundary terms at infinity vanish thanks to the rapid decay  of $\log\det\nolimits_F T^\sharp(\mu)$ proved in  the third assertion of
Lemma~\ref{Lem3.1}. As a result,  with the branch of the logarithm specified in that assertion,
we have
\begin{equation}\label{eq4.9}
 \zeta_{\mathrm{rel}}^\sharp(s;\lambda)
 =
 -\frac{s}{2\pi i}\int_{\mathcal C_\lambda}
 (\lambda-\mu)^{-s-1}
 \log\det\nolimits_F T^\sharp(\mu)\,d\mu.
\end{equation}
The same assertion also ensures locally uniform convergence of
the integral~\eqref{eq4.9} in $s$. It therefore provides an analytic
continuation through $s=0$.

We next establish the regularity at $s=0$ of the zeta functions
$\zeta_m^\sharp(s;\lambda)$.  For $\Re s>1$, the orthogonal
decompositions of the Dirichlet Laplacians $\Delta_{\partial X_\epsilon,h}$ together with the equality
$m=\widetilde m$ on $X_\epsilon$ imply
\begin{equation}\label{eq4.10}
\zeta_{\partial X_\epsilon,m}^\sharp(s;\lambda)
-\zeta_{\partial X_\epsilon,\widetilde m}^\sharp
 (s;\lambda)
=
 \zeta_{\mathcal U,m}^\sharp(s;\lambda)
-\zeta_{\mathcal U,\widetilde m}^\sharp(s;\lambda),
\end{equation}
where
$$
 \zeta_{\mathcal U,h}^\sharp(s;\lambda)=\sum_j  \zeta_{\mathcal U^j,h}^\sharp(s;\lambda),   \qquad  \zeta_{\mathcal U^j,h}^\sharp(s;\lambda)=\Tr(\Delta_{\mathcal U^j,h}^\sharp+\lambda)^{-s}.
$$
The zeta functions $\zeta_{\mathcal U^j,h}^\sharp$ for the Laplacians in the disks are regular at $s=0$
by Proposition~\ref{Prop2.1}; the addition of the
positive constant $\lambda$ does not affect this regularity.
Consequently, the difference of the  zeta functions in the left hand side of
\eqref{eq4.10} is regular at $s=0$.

The metric $\widetilde m$ is smooth on
$X$, and the restrictions of
$\Delta_{\widetilde m}$ to odd and even functions are unitarily
equivalent to the standard Dirichlet and Neumann Laplacians on the
smooth quotient $(M, \widetilde m)$.  Hence
$\zeta_{\widetilde m}^\sharp(s;\lambda)$ is regular at $s=0$.
Solving~\eqref{eq4.7} for the remaining zeta
function leads to
$$
 \zeta_m^\sharp(s;\lambda)
 =\zeta_{\widetilde m}^\sharp(s;\lambda)
 +\left(\zeta_{\partial X_\epsilon,m}^\sharp(s;\lambda)
 -\zeta_{\partial X_\epsilon,\widetilde m}^\sharp
      (s;\lambda)\right)-\zeta_{\mathrm{rel}}^\sharp(s;\lambda),
$$
where the right hand side is regular at $s=0$.
Thus $\zeta_m^\sharp(s;\lambda)$ is regular at $s=0$, and  the zeta regularized determinants in~\eqref{eq4.1} are well-defined. 

Differentiating the contour representation~\eqref{eq4.9} at $s=0$, we obtain
$$
 \partial_s\zeta_{\mathrm{rel}}^\sharp(0;\lambda)
 =
 \frac1{2\pi i}\int_{\mathcal C_\lambda}
 \frac{\log\det\nolimits_F T^\sharp(\mu)}{\mu-\lambda}\,d\mu.
$$
Truncate $\mathcal C_\lambda$ at $\Re\mu=L$ and close it clockwise by the vertical segment. The integral over this segment tends to zero as $L\to+\infty$, see  Lemma~\ref{Lem3.1}. By Cauchy's formula we get \begin{equation}\label{eq4.11}
 \partial_s\zeta_{\mathrm{rel}}^\sharp(0;\lambda)
 =-\log\det\nolimits_F T^\sharp(\lambda).
\end{equation}

%=================

%=================
By the definition of the zeta regularized determinant and the
localization of the zeta-function difference
in~\eqref{eq4.10},  we have  
$$\log \frac{\Det(\Delta_{\widetilde m}^\sharp+\lambda)}
      {\Det(\Delta_m^\sharp+\lambda)}
 \frac{\Det(\Delta_{\mathcal U,m}^\sharp+\lambda)}
      {\Det(\Delta_{\mathcal U,\widetilde m}^\sharp+\lambda)}=
 -\partial_s\zeta_{\mathrm{rel}}^\sharp(0;\lambda).
$$
This together with~\eqref{eq4.11} completes the proof.
 \end{proof}

%-----------=======================----------------------

\section{Determinant comparison under conical smoothing}\label{Sec5}
In this section we complete the proof of Theorem~\ref{Thm1.1}.
To this end, we pass to the limit $\lambda\to0^+$ in the relative BFK
identity of Proposition~\ref{Prop4.1}.
This justifies the determinant comparison identities stated in Proposition~\ref{Prop5.1} below. These identities, together with
the local formula of Proposition~\ref{Prop2.1}, provide
the main ingredients of the proof.
We use the same notation as in Section~\ref{Sec4}.

\begin{proposition}\label{Prop5.1}
The following comparison identities hold:
\begin{equation}\label{eq5.1}
 \frac{\Det\Delta_m^D}
      {\Det\Delta_{\widetilde m}^D}
 =
 \frac{\Det\Delta_{\mathcal U,m}^D}
      {\Det\Delta_{\mathcal U,\widetilde m}^D},
 \qquad
 \frac{\Det\Delta_m^N/\Area(X,m)}
      {\Det\Delta_{\widetilde m}^N/\Area(X,\widetilde m)}
 =
 \frac{\Det\Delta_{\mathcal U,m}^N}
      {\Det\Delta_{\mathcal U,\widetilde m}^N}.
\end{equation}
Here the determinants are zeta regularized in the standard way. 
In particular, the zero eigenvalue of
$\Delta_h^N$ is omitted.
\end{proposition}

The proof  of Proposition~\ref{Prop5.1} essentially reduces to evaluating, as
$\lambda\to0^+$, the limits of the Fredholm determinants of
$T^D(\lambda)$ and $T^N(\lambda)$, the restrictions of the
relative Neumann jump operator to odd and even functions. For odd functions, the restrictions of the Neumann jump operators $R_h^D(\lambda)$, $h\in\{m,\widetilde m\}$, 
are invertible at $\lambda=0$, and the Fredholm determinant of the
relative Neumann jump operator $T^D(\lambda)$  tends to one.  For even functions, the restrictions $R_h^N(0)$ are not invertible and 
the constant zero mode produces an additional ratio of areas in the limit of $\det\nolimits_FT^N(\lambda)$. These assertions are proved in
Lemmas~\ref{Lem5.2} and~\ref{Lem5.3} below.
The proof of Proposition~\ref{Prop5.1} follows the two lemmas.
We then complete the proof of Theorem~\ref{Thm1.1}.

We need some preliminaries before formulating the first lemma. 
For $\varphi\in H^1(\partial X_\epsilon)$, let
$u_h(0)$ denote the unique weak solution of the boundary value problem~\eqref{eq3.1} 
with $\lambda=0$.
The weak harmonic equation in dimension two is conformally invariant,
and hence
$u_m(0)=u_{\widetilde m}(0)$.  Since the smoothed metric $\widetilde m$  coincides with the conical metric $m$ near
$\partial X_\epsilon$, the corresponding normal
derivatives are also equal.  Consequently, for $\lambda=0$ the Neumann jump operators for $h=\widetilde m$ and $h=m$  coincide as unbounded operators in $ L^2(\partial X_\epsilon)$ with domain $H^1(\partial X_\epsilon)$, 
\begin{equation}\label{eq5.2}
 R_{\widetilde m}(0)=R_m(0):
 H^1(\partial X_\epsilon)
 \longrightarrow
 L^2(\partial X_\epsilon);
\end{equation}
cf.~ \cite[Proof of Lemma~2.10]{Ref11}. 

Moreover, for each metric $h\in\{m,\widetilde m\}$, reflection equivariance reads
$ R_h(0)\Pi^\sharp=\Pi^\sharp R_h(0)$ with $\sharp\in\{D,N\}$. Equivalently, $R_h(0)$ preserves the subspaces of odd and even  functions on $\partial X_\epsilon$.
Restriction of \eqref{eq5.2} to
$\Pi^\sharp L^2(\partial X_\epsilon,ds)$ yields 
$$R_{\widetilde m}^\sharp(0)=R_m^\sharp(0),\qquad\sharp\in\{D,N\}.$$

%------------------------------
\begin{lemma}\label{Lem5.2}
For $h\in\{m,\widetilde m\}$, the restriction of the Neumann jump $R_h(0)$ to the
odd functions,
$$
 R_h^D(0):
 \Pi^D H^1(\partial X_\epsilon)
 \longrightarrow
 \Pi^D L^2(\partial X_\epsilon),
$$
is invertible.  Moreover, the Fredholm determinant of the relative Neumann jump $T(\lambda)$ restricted to the odd functions satisfies
\begin{equation}\label{eq5.3}
 \lim_{\lambda\to0}\det\nolimits_FT^D(\lambda)=1.
\end{equation}
\end{lemma}

\begin{proof}
Suppose that
$\varphi\in\Pi^D H^1(\partial X_\epsilon)$ and
$R_h^D(0)\varphi=0$.  The weak solution $u_h(0)$ of the Dirichlet problem~\eqref{eq3.1} has matching normal derivatives across
$\partial X_\epsilon$ and therefore extends to a
global harmonic function on $X$ that is  in the
Friedrichs domain of $\Delta_h$.  This function is constant and,
being odd with respect to $\iota$, vanishes.  Thus the selfadjoint operator  $R_h^D(0)$ in $ \Pi^D L^2(\partial X_\epsilon)$ with domain $\Pi^D H^1(\partial X_\epsilon)$ has a  compact resolvent and trivial kernel. Hence $R_h^D(0)$ is invertible.

By assertion~1 of Lemma~\ref{Lem3.1}, near $\lambda=0$ the operators $\lambda\mapsto R_m^D(\lambda)$ form an analytic  family of type A in the sense of Kato. This together with invertibility of $R_h^D(0)$ implies  that the family
 $$
\lambda\longmapsto  \bigl(R_m^D(\lambda)\bigr)^{-1}\in
\mathcal B\left( \Pi^D L^2(\partial X_\epsilon),\Pi^D H^1(\partial X_\epsilon)\right)
$$
is also analytic near zero~\cite{Ref20}.  The same assertion of that lemma shows that
$$
 \lambda\longmapsto \left(R_{\widetilde m}^D(\lambda)-R_m^D(\lambda)\right)\in\mathcal B\left(
 \Pi^D H^1(\partial X_\epsilon)
,
 \Pi^D H^2(\partial X_\epsilon)\right)
$$
is analytic.  Consequently,
$$
\lambda\mapsto T^D(\lambda)-I
 =
 \bigl(R_{\widetilde m}^D(\lambda)-R_m^D(\lambda)\bigr)
 \bigl(R_m^D(\lambda)\bigr)^{-1}\in
\mathcal B\left( \Pi^D L^2(\partial X_\epsilon),\Pi^D H^2(\partial X_\epsilon)\right)
$$
is analytic.  Since the embedding 
$H^2(\partial X_\epsilon)\hookrightarrow
L^2(\partial X_\epsilon)$ is trace class, this is
an analytic $\mathfrak S_1$-valued family.  Thanks to $R_{\widetilde m}^D(0)=R_m^D(0)$
we have  $T^D(0)=I$, and hence
$$
 \|T^D(\lambda)-I\|_{\mathfrak S_1}\longrightarrow0
 \qquad\text{as }\lambda\to0.
$$
Continuity of the Fredholm determinant in the trace norm proves~\eqref{eq5.3}.
\end{proof}

%----------------------------------------

\begin{lemma}\label{Lem5.3}
Regard $R_h^N(\lambda)$ as an operator in
$\Pi^N L^2(\partial X_\epsilon)$, with domain
$\Pi^N H^1(\partial X_\epsilon)$.  It is selfadjoint for real $\lambda$.  Its kernel at
$\lambda=0$ consists of the constant functions.  If
$\mu_h^N(\lambda)$ denotes the unique eigenvalue converging to zero,
then
\begin{equation}\label{eq5.4}
 \mu_h^N(\lambda)
 =
 \lambda\frac{\Area(X,h)}
 {\Length(\partial X_\epsilon)}
 +O(\lambda^2),\quad \lambda\to 0.
\end{equation}
Moreover,  the Fredholm determinant of the relative Neumann jump $T(\lambda)$ restricted to the even functions satisfies
\begin{equation}\label{eq5.5}
 \lim_{\lambda\to0}\det\nolimits_FT^N(\lambda)
 =
 \frac{\Area(X,\widetilde m)}
 {\Area(X,m)}.
\end{equation}
\end{lemma}

\begin{proof} The same argument as in the beginning of the proof of Lemma~\ref{Lem5.2} implies that the kernel of 
 $R_h^N(0)$ is spanned by $1_{\partial X_\epsilon}$.

By assertion~(1) of Lemma~\ref{Lem3.1}, the family
$R_h^N(\lambda)$ is analytic near zero.
By  analytic perturbation theory~\cite{Ref20},  the simple zero eigenvalue
extends to  an analytic eigenvalue $\lambda\mapsto \mu_h^N(\lambda)$ of $R_h^N(\lambda)$ with  an analytic eigenfunction $\lambda\mapsto e_h(\lambda)$ satisfying$$
 R_h^N(\lambda)e_h(\lambda)
 =\mu_h^N(\lambda)e_h(\lambda),
 \qquad e_h(0)=1_{\partial X_\epsilon}.
$$
Differentiating this identity at $\lambda=0$ and using
$\mu_h^N(0)=0$, we obtain
$$
 \partial_\lambda R_h^N(0)  1_{\partial X_\epsilon}
 +R_h^N(0)\partial_\lambda e_h(0)
 = \partial_\lambda\mu_h^N(0)1_{\partial X_\epsilon}.
$$
Next, we take the inner product in $L^2(\partial X_\epsilon)$ of both sides of this equality with $1_{\partial X_\epsilon}$. The second term on the left
vanishes because $R_h^N(0)$ is selfadjoint and $R_h^N(0)1_{\partial X_\epsilon}=0$.
The invertibility of
$\Delta_{\partial X_\epsilon,h}$ also implies
that $ \partial_\lambda R_h^N(0)
 =
 \left.
 \mathcal P_h^*(0)\mathcal P_h(0)
 \right|_{\Pi^N L^2(\partial X_\epsilon)}$   as  both sides
of~\eqref{eq4.4} are analytic near zero. 
Consequently, we have
$$
 \Area(X,h)\!=\bigl\|\mathcal P_h(0) 1_{\partial X_\epsilon}\bigr\|_{L^2(X,h)}^2=\bigl\langle
       \partial_\lambda R_h^N(0) 1_{\partial X_\epsilon},  1_{\partial X_\epsilon}\bigr\rangle_{L^2(\partial X_\epsilon)}=\partial_\lambda\mu_h^N(0)  \Length(\partial X_\epsilon),
$$
where  the Poisson extension $\mathcal P_h(0)1_{\partial X_\epsilon} $ of $1_{\partial X_\epsilon}$
is the global constant $1$ on $(X,h)$. This together with analyticity of $\lambda\mapsto \mu_h^N(\lambda)$ proves \eqref{eq5.4}.

%%=========

 It remains to evaluate the limit of the Fredholm determinant in~\eqref{eq5.5}.  Let $\mathsf P^N_h(\lambda)$ stand for  the Riesz projection onto the kernel  $\ker\left( R_h^N(\lambda)-\mu_h^N(\lambda)\right)$.  
Then 
$$
\mathsf P^N:=\mathsf P^N_m(0)=\mathsf P^N_{\widetilde m}(0)
$$
 is the orthogonal projection in $
\Pi^N L^2(\partial X_\epsilon)$ onto $\ker R_h^N(0)
$.

Clearly, the operator $R_m^N(0)+\mathsf P^N$ is invertible. As a consequence, in a small neighborhood of $\lambda=0$, we have
$$
 \bigl(R_m^N(\lambda)\bigr)^{-1}
 =
 \frac{1}{\mu_m^N(\lambda)} \mathsf P^N_m(\lambda)
 +
 \bigl(R_m^N(\lambda)+\mathsf P^N_m(\lambda)\bigr)^{-1}
 \bigl(I-\mathsf P^N_m(\lambda)\bigr).
$$
The expansion
$\mathsf P^N_m(\lambda)=\mathsf P^N+O(\lambda)$ together with the expansion~\eqref{eq5.4} for the eigenvalue $ \mu_h^N(\lambda)$ implies
$$
 \bigl(R_m^N(\lambda)\bigr)^{-1}
 =
 \frac{ {\Length(\partial X_\epsilon)}}{\Area(X,m)}\frac 1 \lambda\mathsf  P^N+G_m(\lambda),
$$
where $G_m(\lambda)$ is an operator analytic near zero from
$\Pi^N L^2(\partial X_\epsilon)$ to
$\Pi^N H^1(\partial X_\epsilon)$.

Consequently, for sufficiently small nonzero $\lambda$ the restriction $T^N(\lambda)$ of the relative Neumann jump operator to the even functions satisfies
\begin{equation}\label{eq5.6}
 T^N(\lambda)-I
 =
 \frac{ {\Length(\partial X_\epsilon)}}{\Area(X,m)}\frac 1 \lambda\bigl(R_{\widetilde m}^N(\lambda)
             -R_m^N(\lambda)\bigr)\mathsf P^N
             +
 \bigl(R_{\widetilde m}^N(\lambda)-R_m^N(\lambda)\bigr)
 G_m(\lambda).
\end{equation}
Since $R_{\widetilde m}^N(0)=R_m^N(0)$,
  the second term on the right  is $O(\lambda)$ in the operator norm of
$\mathcal B\left(\Pi^N L^2(\partial X_\epsilon),\Pi^N H^2(\partial X_\epsilon)\right)$; see  Lemma~\ref{Lem3.1}. Hence, it is  also $O(\lambda)$ in the trace norm in   $\Pi^N L^2(\partial X_\epsilon)$,  as the embedding $H^2(\partial X_\epsilon)\hookrightarrow
L^2(\partial X_\epsilon)$ is trace class.
Thus
$$
 T^N(\lambda)-I
 \to   \frac{ {\Length(\partial X_\epsilon)}}{\Area(X,m)}\bigl(\partial_\lambda R_{\widetilde m}^N(0)
             -\partial_\lambda R_m^N(0)\bigr)\mathsf P^N
 \quad\text{in}\quad \mathfrak S_1(\Pi^N L^2(\partial X_\epsilon)) \quad\text{as}\quad \lambda\to 0.
$$
The representation~\eqref{eq5.6} also shows that $T^N(\lambda)-I$ extends analytically to $\lambda=0$ in the trace norm.
Finally, thanks to the continuity of the Fredholm determinant,  we  obtain
$$
 \begin{aligned}
 \lim_{\lambda\to0}\det\nolimits_F T^N(\lambda)
 &=
 1+\frac{
 \left\langle
 \bigl(\partial_\lambda R_{\widetilde m}^N(0)
       -\partial_\lambda R_m^N(0)\bigr) 1_{\partial X_\epsilon}, 1_{\partial X_\epsilon}
 \right\rangle}{\Area(X,m)}\\
 &=1+
 \frac{\Area(X,\widetilde m)-\Area(X,m)}{\Area(X,m)}
 =
 \frac{\Area(X,\widetilde m)}
      {\Area(X,m)}.
 \end{aligned}
$$
This completes the proof.
\end{proof}

\begin{proof}[Proof of Proposition~\ref{Prop5.1}]

%=======

We first show that $\Det\Delta_m^D$ and $\Det\Delta_m^N$ are well defined, since Section~\ref{Sec4} establishes zeta regularity for $\lambda>0$ only.

The proofs of Lemmas~\ref{Lem5.2} and~\ref{Lem5.3}
show that the branch of $\log\det_F T^\sharp$ specified in
Lemma~\ref{Lem3.1} extends analytically to a neighborhood of zero.
In the representation of the relative zeta
function~\eqref{eq4.9}, we deform the contour
$\mathcal C_\lambda$ within this neighborhood to a contour
$\mathcal C$ passing to the left of zero.
For all sufficiently small $\lambda>0$, the contour
$\mathcal C$ can be chosen to be $\lambda$-independent.
Since the integrand is analytic between the two contours, we conclude from Cauchy's theorem that
$$
\zeta_{\mathrm{rel}}^\sharp(s;\lambda)
=
-\frac{s}{2\pi i}\int_{\mathcal C}
(\lambda-\mu)^{-s-1}
\log\det\nolimits_F T^\sharp(\mu)\,d\mu.
$$

Let $r=\operatorname{dist}\{0,\mathcal C\}>0$. On the contour $\mathcal C$, we have
$$
|\lambda-\mu|\ge  r/2,
\qquad \lambda\in[0, r/2].
$$
Thus the integrand converges as $\lambda\to0^+$, uniformly for $\mu$ on any finite part of $\mathcal C$ and $s$ in any compact subset of $\mathbb C$. On the two half-lines of $\mathcal C$, Lemma~\ref{Lem3.1} supplies an integrable majorant, uniform for $\lambda\in[0, r/2]$ and for $s$ in any compact subset of $\mathbb C$. We may therefore pass to the limit under the integral:
$$
\zeta_{\mathrm{rel}}^\sharp(s;\lambda)
\longrightarrow
\zeta_{\mathrm{rel}}^\sharp(s;0)
:=
-\frac{s}{2\pi i}\int_{\mathcal C}
(-\mu)^{-s-1}
\log\det\nolimits_F T^\sharp(\mu)\,d\mu,
$$
locally uniformly in $s$. In particular, the limit is analytic near $s=0$, and
$$
\partial_s\zeta_{\mathrm{rel}}^\sharp(0;\lambda)
\longrightarrow
\partial_s\zeta_{\mathrm{rel}}^\sharp(0;0).
$$

%=======

In the relative zeta function   $\zeta_{\mathrm{rel}}^N(s;\lambda)$ the contributions $\lambda^{-s}$ of the constant modes of $\Delta_h^N$ cancel, cf.~\eqref{eq4.7}.
Thanks to~\eqref{eq4.10}, we have
$$
\begin{aligned}
 {}&\zeta_m^D(s;\lambda)
 =\zeta_{\widetilde m}^D(s;\lambda)
   +\zeta_{\mathcal U,m}^D(s;\lambda)-\zeta_{\mathcal U,\widetilde m}^D(s;\lambda)
   -\zeta_{\mathrm{rel}}^D(s;\lambda),\\
 {} & \zeta_m^N(s;\lambda)-\lambda^{-s}=\zeta_{\widetilde m}^N(s;\lambda)-\lambda^{-s}
   +\zeta_{\mathcal U,m}^N(s;\lambda)-\zeta_{\mathcal U,\widetilde m}^N(s;\lambda)
   -\zeta_{\mathrm{rel}}^N(s;\lambda).
\end{aligned}
$$
For $\Re s>1$, the two left hand sides converge as $\lambda\to0^+$ to the spectral zeta functions of $\Delta_m^D$ and $\Delta_m^N$, respectively.
The right hand sides converge locally uniformly in $s$ near $s=0$
to analytic limits, which therefore give the analytic continuations of
these spectral zeta functions to $s=0$. The derivatives at zero
converge as well. Thus the determinants $\Det\Delta_m^\sharp$ are well defined.

As  $\lambda\to0^+$, the determinants satisfy 
$$
\begin{aligned}
 \Det(\Delta_h^D+\lambda)
 &=\Det\Delta_h^D\,(1+o(1)),\\
 \Det(\Delta_h^N+\lambda)
 &=\lambda\Det\Delta_h^N\,(1+o(1)),\\
 \Det(\Delta_{\mathcal U,h}^\sharp+\lambda)
 &=\Det\Delta_{\mathcal U,h}^\sharp\,(1+o(1)),
\end{aligned}
$$
where $h\in\{m,\widetilde m\}$ and  $\sharp\in\{D,N\}$. Passing to the limit in the relative BFK identity from Proposition~\ref{Prop4.1} and
using Lemma~\ref{Lem5.2} implies the first identity in~\eqref{eq5.1}.  The same limit with
Lemma~\ref{Lem5.3} proves the second identity in~\eqref{eq5.1}. \end{proof}

Let $h$ be either the  metric $m$ or a smooth reflection invariant conformal
metric on $X$. We shall use the same notation with subscript $h$
for the corresponding Laplacians. The Laplacian $\Delta_h$ commutes
with the orthogonal projectors onto odd and even functions in
$L^2(X,h)$. Its restrictions to these subspaces are unitarily
equivalent to $\Delta_h^D$ and $\Delta_h^N$, respectively. Therefore 
\begin{equation}\label{eq5.7}
\Det\Delta_h=\Det\Delta_h^D\Det\Delta_h^N.
\end{equation}

\begin{lemma}[After W. L\"uck and S.\,M. Vishik]\label{Lem5.4}
Let $h$ be a smooth conformal metric on $M$ which extends by
reflection to a smooth metric on $X$. Then
$$
\frac{\Det\Delta_h^N}
 {\Area(M,h)\Det\Delta_h^D}
 =\frac{\|v\|_{\tau_0(M)}^2}
 {\Area(M,m)\|v\|_{\det H^\bullet(M)}^2}.
$$
\end{lemma}

\begin{proof}
In this proof, we write $\|\cdot\|_{\det H^\bullet(M,h)}$ for the
norm induced by the metric $h$, and $T_0(M,h)$ for the
corresponding analytic torsion norm.
Choose a smooth conformal metric $m_0$ on $M$ which is a
product near the boundary. Such a metric is obtained by taking
$dt^2+ds^2$ in conformal collars of the boundary components and
extending it smoothly over $M$ within its conformal class.
For the product metric $m_0$, Vishik's theorem~\cite[Theorem~1.4]{Ref30}
implies
$$
 \|v\|_{T_0(M,m_0)}^2=\|v\|_{\tau_0(M)}^2
$$
for any nonzero $v\in\det H^\bullet(M)$.
%----------------

Let $\zeta_q$ denote the spectral zeta function of the de Rham
Laplacian on $q$-forms with absolute boundary conditions; as usual, the
zero eigenvalues are omitted. The Hodge decomposition and the Hodge
star imply
$$
 \zeta_0=\zeta_{m_0}^N,\qquad
 \zeta_1=\zeta_{m_0}^N+\zeta_{m_0}^D,\qquad
 \zeta_2=\zeta_{m_0}^D,
$$
where $\zeta_{m_0}^D$ and $\zeta_{m_0}^N$ are the spectral zeta
functions of the scalar Dirichlet and Neumann Laplacians on
$(M,m_0)$, respectively.
Indeed, for each positive eigenvalue, the one-form eigenspace is
spanned by differentials of Neumann eigenfunctions and Hodge stars
of differentials of Dirichlet eigenfunctions with that eigenvalue.
The Hodge star also identifies two-form eigenforms satisfying
absolute boundary conditions with Dirichlet eigenfunctions. Thus Vishik's definition~\cite[(1.2)]{Ref30} yields
$$
 \|v\|_{T_0(M,m_0)}^2
 =\|v\|_{\det H^\bullet(M,m_0)}^2
   \exp\left(\sum_{q=0}^2(-1)^q q\zeta_q'(0)\right)
 =\|v\|_{\det H^\bullet(M,m_0)}^2
   \frac{\Det\Delta_{m_0}^N}{\Det\Delta_{m_0}^D}.
$$
As a consequence we have
\begin{equation}\label{eq5.8}
\frac{\Det\Delta_{m_0}^N}
 {\Area(M,m_0)\Det\Delta_{m_0}^D}
 =\frac{\|v\|_{\tau_0(M)}^2}
 {\Area(M,m_0)\|v\|_{\det H^\bullet(M,m_0)}^2}.
\end{equation}

By the factorization~\eqref{eq5.7} and the equality
$\Area(X,m_0)=2\Area(M,m_0)$, the quotient on the left equals
$$2\,\dfrac{\Det\Delta_{m_0}/\Area(X,m_0)}
 {(\Det\Delta_{m_0}^D)^2}.$$
Hence, within the class of conical conformal metrics admissible
in the sense of~\cite[Definition~2.1]{Ref11} and smooth near
$\partial M$, the conformal invariance of the 
left hand side in~\eqref{eq5.8} follows from the Polyakov formula on $X$ and the
Polyakov--Alvarez formula on its two isometric halves;
see~\cite[(1.1) and (1.4)]{Ref11}.  On the right hand side of~\eqref{eq5.8}, the area factor cancels the contribution of
$H^0(M)$ to the squared determinant-line norm. The remaining norm
on $\det H^1(M)$ is conformally invariant. Indeed,  the Hodge star on
one-forms is conformally invariant, while the combinatorial torsion
norm is independent of the metric.

We may therefore replace $m_0$ by $h$ on the left hand side
of~\eqref{eq5.8} and by $m$ on the right hand side. This completes the proof.
\end{proof}

\begin{proof}[Proof of Theorem~\ref{Thm1.1}]
If $\partial M$ has no corners, the result follows directly from
Lemma~\ref{Lem5.4} and the spectral decomposition
into odd and even functions~\eqref{eq5.7}.

If $\partial M$ has  corners, we combine Proposition~\ref{Prop5.1} with
Proposition~\ref{Prop2.1}.  Since on the union $\mathcal U$ of the disks $\mathcal U^j$ we have
$\Delta_{\mathcal U,h}^\sharp=\bigoplus_j
\Delta_{\mathcal U^j,h}^\sharp$, the determinant of  $\Delta_{\mathcal U,h}^\sharp$ is the product of
the determinants of $\Delta_{\mathcal U^j,h}^\sharp$ on the individual disks.  We use $F_j(h)$ as
defined in Section~\ref{Sec2}.
Dividing the second identity in~\eqref{eq5.1} by
the first one and using~\eqref{eq2.2} and the equality
$\Area(X,h)=2\Area(M,h)$, we obtain
\begin{equation}\label{eq5.9}
 \begin{aligned}
 \frac{\Det\Delta_m^N}
      {\Area(M,m)\Det\Delta_m^D}
 &=
 \frac{\Det\Delta_{\widetilde m}^N}
      {\Area(M,\widetilde m)\Det\Delta_{\widetilde m}^D}
 \prod_j\frac{F_j(m)}{F_j(\widetilde m)}
 \\
 &=
 \frac{\Det\Delta_{\widetilde m}^N}
      {\Area(M,\widetilde m)\Det\Delta_{\widetilde m}^D}
 \prod_j\alpha_j^{-1/2}.
 \end{aligned}
\end{equation}

Substituting the identity from Lemma~\ref{Lem5.4} with $h=\widetilde m$ into~\eqref{eq5.9} and multiplying by $\Area(M,m)$, we obtain
$$
 \frac{\Det\Delta_m^N}{\Det\Delta_m^D}
 =\frac{\|v\|_{\tau_0(M)}^2}
 {\|v\|_{\det H^\bullet(M)}^2}
 \prod_j\alpha_j^{-1/2}.
$$
Combining this identity with the factorization~\eqref{eq5.7} and
solving for  the determinants $\Det\Delta_m^D$ and $\Det\Delta_m^N$ proves our main result~\eqref{eq1.1}.
\end{proof}

%------------------------------------------

\section{BFK gluing and the torsion factor}\label{Sec6}
For a smooth reflection invariant metric $h$ on $X$, conformal to $m$, let $\mathcal R_h$ denote 
the Neumann jump operator for the cut of $X$ along $\partial M$.  In this smooth case, the following three expressions are equal and define the same conformal invariant:
$$
\frac{\Det\mathcal R_h}{\Length(\partial M,h)}
=
\frac{\Det\Delta_h/\Area(X,h)}
     {(\Det\Delta_h^D)^2}
=
\frac12\,
\frac{\|v\|_{\tau_0(M)}^2}
     {\Area(M,h)\,
      \|v\|_{\det H^\bullet(M,h)}^2},
$$
where the first equality is the BFK formula~\cite{Ref5} and the second follows from Lemma~\ref{Lem5.4}.  
For the reflection invariant conical metrics admissible in the sense
of~\cite[Definition~2.1]{Ref11} and smooth near $\partial M$,
all three expressions remain well defined and equal;
see~\cite[Proposition~2.8]{Ref11} for the BFK equality. 
%=============

The situation changes when conical singularities of the double $
 X=M\cup_{\partial M}M$ lie on $\partial M$, producing corners on the boundary of each half.
The last expression can be evaluated using 
 any finite-area conical metric in the same conformal class, without changing its value; see the proof of Lemma~\ref{Lem5.4}. The second expression remains well defined and, by
Theorem~\ref{Thm1.1}, equals the third expression multiplied by
$\prod_j\alpha_j^{-1/2}$. It is therefore conformally invariant for
fixed corner angles. For the first expression, the usual
definition of the Neumann jump determinant as the zeta regularized determinant of a  classical elliptic pseudodifferential operator~\cite{Ref27} no longer
applies, since the cut passes through conical points.

%===========

The expression involving the Reidemeister torsion norm can also be written in terms of harmonic one-forms.

Indeed, collapse a triangulation of $M$ to a graph and contract a spanning tree.
Each cancelled cellular pair has coefficient $\pm1$, so these reductions
preserve the Reidemeister norm.
Choose the corresponding loops as the cycles $\gamma_i$, and let
$\omega_j$ be the harmonic one-forms satisfying absolute boundary
conditions and normalized by
$$
\int_{\gamma_i}\omega_j=\delta_{ij}.
$$
Introduce the Gram matrix $G$ with entries
\begin{equation}\label{eq6.1}
G_{ij}=\int_M\omega_i\wedge *\overline{\omega_j}.
\end{equation}
The reduced complex has zero differential
and orthonormal cochains representing $1$ and the cohomology classes $[\omega_j]$. Thus, for  $v=1^{-1}\otimes\bigwedge_j[\omega_j]$,  Vishik's definition~\cite[(1.6)]{Ref30} implies
$\|v\|_{\tau_0(M)}^2=1$. On the other hand, we have 
$$
\|v\|_{\det H^\bullet(M,m)}^2
=\|1\|_{L^2(M,m)}^{-2}\det G
=\frac{\det G}{\Area(M,m)}.
$$
Consequently,
$$
\frac12\,
\frac{\|v\|_{\tau_0(M)}^2}
{\Area(M,m)\,
 \|v\|_{\det H^\bullet(M,m)}^2}
=\frac1{2\det G};
$$
see also~\cite[Proposition~5.4 and (5.2)]{Ref26}.
The empty Gram determinant is understood to be equal to $1$. In particular, if $M$ is a topological disk, then the quotient $ \frac{\|v\|_{\tau_0(M)}^2}
{ \|v\|_{\det H^\bullet(M)}^2}$ in Theorem~\ref{Thm1.1} reduces to the area of the polygon
$(M,m)$. In this case,  the first formula in our main result~\eqref{eq1.1} was announced in conjectural form in~\cite[Section~6]{Ref16}. For circular cuts of the Riemann sphere avoiding the conical points, the value $1/2$ of the normalized Neumann jump determinant was obtained in~\cite[(2.6)]{Ref14}. This computation does not require reflection symmetry of the metric.

The use of the double to compare analytic and combinatorial torsion already appears in~\cite{Ref26,Ref30}.
For a reflection invariant metric $h$ smooth near $\partial M$, one has  $\Det\mathcal R_h=\frac12\Det\mathcal N_h$, where $\mathcal N_h$ is the Dirichlet-to-Neumann operator of either half.  It was proved in~\cite{Ref8} that  $\Det \mathcal N_h/\operatorname{Length}(\partial M,h)=1$ for smooth simply connected planar domains. Conformal invariance of this quotient for compact surfaces with smooth
boundary and its expression through Ruelle and Selberg zeta functions were established in~\cite{Ref9}. For the same quotient, formulas in terms of periods were obtained in~\cite{Ref23,Ref21} for surfaces with connected
smooth boundary and for multiply connected
planar domains with smooth boundary, respectively. By conformal invariance, these formulas follow from
Lemma~\ref{Lem5.4} and the BFK formula after
expressing the Gram matrix~\eqref{eq6.1} through the periods
of the double.

%--------------------------

\section{Examples and applications}\label{Sec7}

\noindent\textbf{The Aurell--Salomonson formula for polygons.}
As a first example, we deduce the Aurell--Salomonson
formula~\cite{Ref3,Ref4} for Euclidean polygons.
To this end, we apply~\cite[ Proposition~3.3]{Ref11} to the metric
double and then use Theorem~\ref{Thm1.1}.

Let $\overline{\mathbb H}=\{z\in\mathbb C:\operatorname{Im}z\geq0\}\cup\{\infty\}$
denote the extended closed upper half-plane.
In a suitably normalized Schwarz--Christoffel coordinate, a bounded
simply connected Euclidean polygon is isometric to
$(\overline{\mathbb H},m)$, where
\begin{equation}\label{eq7.1}
 m=\prod_{j=1}^N|z-p_j|^{2\alpha_j-2}|dz|^2,\quad p_j\in\mathbb R,
\end{equation}
with $N\geq3$,  $ \alpha_j\in(0,2)$, and 
 $\sum_{j=1}^N\alpha_j=N-2$.
The  points $p_j\in\mathbb R$ are distinct and correspond to the
vertices with interior angles $\pi\alpha_j$.

For $\alpha>0$, define
\begin{equation}\label{eq7.2}
 C(\alpha):=2\zeta_B'(0;\alpha,1,1)-2\zeta_R'(-1)
 -\frac{(\alpha-1)^2}{6\alpha}\log2
-\frac{\alpha-1}{12}
 +\frac12\log\alpha.
\end{equation}
Here $\zeta_B(s;a,b,x)$ and $\zeta_R(s)$ are the Barnes double
and Riemann zeta functions, respectively, and the primes denote
differentiation with respect to $s$.

The metric double is the Riemann sphere $X=\overline{\mathbb C}$
with the metric~\eqref{eq7.1}, cone angles
$2\pi\alpha_j$, and
$\Area(X,m)=2\Area(\mathbb H,m)$.
The Aurell--Salomonson formula~\cite[(50)]{Ref4} for
polyhedral spheres was first rigorously derived as a byproduct
in~\cite[Proposition~3.3]{Ref11}. Applied to the double, it takes the form
\begin{equation}\label{eq7.3}
\begin{aligned}
 \log\frac{\Det\Delta_m}{\Area(X,m)}
={}&\frac16\sum_{j=1}^N
       \sum_{\substack{i=1\\ i\ne j}}^N
       \frac{(\alpha_i-1)(\alpha_j-1)}{\alpha_j}
       \log|p_i-p_j|
 -\sum_{j=1}^N C(\alpha_j)\\
 &-\frac43\log2-4\zeta_R'(-1)
 +\frac16-\log\pi.
 \end{aligned}
\end{equation}
For a topological disk, the first identity in Theorem~\ref{Thm1.1}
reads
\begin{equation}\label{eq7.4}
 \log\Det\Delta_m^D
 =\frac12\log\frac{\Det\Delta_m}{\Area(X,m)}
  +\frac12\log2
  +\frac14\sum_{j=1}^N\log\alpha_j.
\end{equation}
Substituting~\eqref{eq7.3} into
\eqref{eq7.4}, we obtain the Aurell--Salomonson
formula~\cite[(32)]{Ref4} for all Euclidean polygons:
\begin{equation}\label{eq7.5}
\begin{aligned}
 \log\Det\Delta_m^D
={}&\frac1{12}\sum_{j=1}^N
       \sum_{\substack{i=1\\ i\ne j}}^N
       \frac{(\alpha_i-1)(\alpha_j-1)}{\alpha_j}
       \log|p_i-p_j|
 -\frac12\sum_{j=1}^N C(\alpha_j)\\
 &+\frac14\sum_{j=1}^N\log\alpha_j
 -2\zeta_R'(-1)
 +\frac1{12}-\frac16\log2-\frac12\log\pi .
\end{aligned}
\end{equation}

For comparison with the original notation and expressions in~\cite{Ref3,Ref4},
see~\cite[Remark~3.4]{Ref11}.
To the best of our knowledge, this is the first rigorous proof of the
Aurell--Salomonson formula for Euclidean polygons.

\medskip
\noindent\textbf{Geodesic triangles.}
Let $(T,m)$ be a spherical or hyperbolic geodesic triangle of constant
curvature $K\ne0$, with interior angles $\pi\alpha_j$, where
$0<\alpha_j<1$. 
By Gauss--Bonnet, its metric double has area
$$
S=2\Area(T,m)
=\frac{2\pi}{K}(\alpha_1+\alpha_2+\alpha_3-1).
$$
The determinant $\Det\Delta_m$ on the double is given explicitly
by~\cite[Corollary~1.3]{Ref14}, with $\beta_j=\alpha_j-1$ and this
value of $S$. Theorem~\ref{Thm1.1} therefore leads to
$$
\begin{aligned}
\Det\Delta_m^D
&=\left(
\frac{K\,\Det\Delta_m}
{\pi(\alpha_1+\alpha_2+\alpha_3-1)}
\right)^{1/2}
\prod_{j=1}^3\alpha_j^{1/4},\\
\Det\Delta_m^N
&=\left(
\frac{\pi(\alpha_1+\alpha_2+\alpha_3-1)}{K}
\,\Det\Delta_m
\right)^{1/2}
\prod_{j=1}^3\alpha_j^{-1/4}.
\end{aligned}
$$

Among geodesic triangles with interior angles less than $\pi$, fixed
curvature $K$, and fixed area $A$, the equilateral triangle is a stationary
point of both determinants, by Theorem~\ref{Thm1.1}
and~\cite[Proposition~5.1]{Ref14}. For each fixed $K$, it is a strict local
minimum of both determinants when $A$ is sufficiently small. Among
Euclidean triangles of unit area, however, it is a strict local maximum
of the Dirichlet determinant and a strict local minimum of the Neumann
determinant.

Indeed, let $a=(1+KA/\pi)/3$, and let $q(a,S)\sum_j(\delta\alpha_j)^2$ be the second
variation of the logarithm of the determinant on the double of the equilateral triangle.
Here $S$ is the area of the double, and the variations $\delta\alpha_j$ satisfy $\sum_j\delta\alpha_j=0$ because curvature and area are fixed.
Theorem~\ref{Thm1.1} implies
\begin{equation}\label{eq7.6}
\delta^2\log\det\Delta_m^{\sharp}
=
\left(\frac12q(a,2A)\mp\frac1{4a^2}\right)
\sum_{j=1}^3(\delta\alpha_j)^2, \quad \sharp\in\{D,N\},
\end{equation}
where the upper sign corresponds to Dirichlet conditions. The scaling
formula and second derivative expressions
in~\cite[(2.3)--(2.4) and Proposition~5.2]{Ref14} give
$$q(a,2A)=-\log(2A)/(6a^3)+O(1),\quad A\to0,$$ with $K$ fixed.
Since $a\to1/3$, both second variations are positive definite for sufficiently
small $A$, so the equilateral triangle is a strict local minimum of both determinants.

At $K=0$ and $A=1$, substitution into~\cite[Proposition~5.2]{Ref14},
using the remainder estimate~\cite[(A.3)]{Ref12} with $N=2$, gives
$-1/4<q(1/3,2)<0$. The second variation~\eqref{eq7.6} is therefore negative definite in the Dirichlet case and positive definite in the Neumann case.

Let us also notice that for triangles of fixed curvature $K\le0$ and area $A>0$, both Dirichlet and Neumann determinants tend to $+\infty$ whenever any interior angle tends to zero. This follows from Theorem~\ref{Thm1.1} and the explicit formula for the determinant on the double of the triangle~\cite[Corollary~1.3]{Ref14}, together with the Barnes zeta asymptotic~\cite[(A.1)]{Ref12}. Consequently, each determinant attains an absolute minimum on a nondegenerate triangle. In particular, the absolute minimum of the Dirichlet determinant among Euclidean triangles of unit area is attained on a non-equilateral triangle.

\medskip
\noindent\textbf{Further examples.}
Closed explicit formulae are available for the determinants of the regular dihedra in~\cite[Examples~5.2 and~5.7]{Ref15}, which are the metric
doubles of regular geodesic polygons. More generally, the closed explicit determinant
formula in~\cite[Theorem~2.1]{Ref15} can be used in
Theorem~\ref{Thm1.1} whenever the triangulated sphere is the metric
double of a geodesic polygon. 

For fixed curvature $K$ and area $A$, the regular
Euclidean and hyperbolic $N$-gons, $N\ge3$, turn out to be  stationary points of both
the Dirichlet and Neumann determinants. Indeed, the corresponding doubles are
stationary by~\cite[Theorem~4.6(2)]{Ref15}. The torsion quotient in
Theorem~\ref{Thm1.1} equals the fixed area $A$, while the Gauss--Bonnet
theorem fixes $\sum_{j=1}^N\alpha_j=N-2+KA/\pi$. At equal angles, the first
variation of $\sum_{j=1}^N\log\alpha_j$ therefore vanishes. The additional
factors are thus stationary, and Theorem~\ref{Thm1.1} transfers stationarity
from the double to both the  Dirichlet and Neumann determinants.

For doubles of genus greater than one, one can apply the determinant
formula~\cite[Proposition~3.9]{Ref11}, with the universal constant
$\kappa_0$ evaluated in~\cite[Remark~7.2, (7.5)]{Ref17}.
Theorem~\ref{Thm1.1} of this paper then allows us to find the Dirichlet
and Neumann determinants without undetermined multiplicative constants.

%
%
%
%Further examples are provided by the explicitly triangulated
%constant curvature conical spheres
%in~\cite[Theorem~2.1 and (2.7)--(2.11)]{Ref15}.
%Whenever such a sphere is the metric double of a geodesic polygon,
%the corresponding Dirichlet and Neumann determinants can be found
%from Theorem~\ref{Thm1.1}.
%
%For fixed curvature and area, the equilateral triangle is a stationary
%point of both the Dirichlet and Neumann determinants. For every
%$N\ge3$, the same holds for regular Euclidean and hyperbolic
%$N$-gons. The corresponding doubles are stationary
%by~\cite[Proposition~5.1]{Ref14}
%and~\cite[Theorem~4.6(2)]{Ref15}, respectively. Indeed, for an
%$N$-gon of curvature $K$, area $A$, and interior angles
%$\pi\alpha_j$, the Gauss--Bonnet theorem implies
%$\sum_{j=1}^N\alpha_j=N-2+{KA}/{\pi}$.
%At equal angles, the variation of the sum $\sum_{j=1}^N\log\alpha_j$ is zero. 
%Since the torsion quotient in Theorem~\ref{Thm1.1} equals the
%fixed area $A$, that theorem transfers stationarity from the double
%to each of the two determinants.
%
%For doubles of genus greater than one, one can apply the determinant
%formula~\cite[Proposition~3.9]{Ref11}, with the universal constant
%$\kappa_0$ evaluated in~\cite[Remark~7.2, (7.5)]{Ref17}.
%Theorem~\ref{Thm1.1} of this paper then allows us to find the Dirichlet and Neumann
%determinants without undetermined multiplicative constants.
%====================================================
\section*{Acknowledgments}
The author thanks Werner M\"uller for a stimulating question on extremal properties of Laplacian determinants on hyperbolic polygons and Paul Wiegmann for helpful discussions.

\end{document}